\documentclass[10pt,leqno]{article}
\usepackage[OT2,OT1]{fontenc}
\usepackage[russian,english]{babel}
\usepackage{a4}
\usepackage{amssymb,amsmath,amsthm}
\usepackage[matrix,arrow]{xy}

\makeatletter

\advance\smallskipamount by 2pt
\advance\medskipamount   by 3pt
\advance\bigskipamount   by 3pt
\normalbaselines
\renewcommand\part{\par\vskip 0pt plus 12\baselineskip \relax
                   \goodbreak
                   \vskip 0pt plus -12\baselineskip \relax
                   \@startsection {part}{0}{\z@}%
                                   {-5ex \@plus -1.5ex \@minus -.5ex}%
                                   {2.3ex \@plus.2ex}%
                                   {\normalfont\Large\bfseries\rightskip=0pt plus 2em\relax}}%
\newcommand*\partmark[1]{}

\renewcommand\section{\@startsection {section}{1}{\z@}%
                                   {-3.5ex \@plus -1ex \@minus -.2ex}%
                                   {2.3ex \@plus.2ex}%
                                   {\normalfont\large\bfseries}}
\renewcommand\subsection{\@startsection{subsection}{2}{\z@}%
                                     {-3.25ex\@plus -1ex \@minus -.2ex}%
                                     {1.5ex \@plus .2ex}%
                                     {\normalfont\normalsize\bfseries}}

\newtheoremstyle{styleDef}{10pt}{10pt}{\upshape}{}{\scshape}{.}{7pt}{\thmnumber{{\upshape\bfseries#2.}\ }\thmname{#1}\thmnote{ #3}}
\newtheoremstyle{styleThm}{10pt}{10pt}{\slshape}{}{\scshape}{.}{7pt}{\thmnumber{{\upshape\bfseries#2.}\ }\thmname{#1}\thmnote{ #3}}
\newtheoremstyle{styleRem}{10pt}{10pt}{\upshape}{}{\itshape}{.}{7pt}{\thmnumber{{\upshape\bfseries#2.}\ }\thmname{#1}\thmnote{ #3}}

\theoremstyle{styleDef}
  \newtheorem{Definition}{Definition}

\theoremstyle{styleThm}
  \newtheorem{Theorem}[Definition]{Theorem}
  \newtheorem*{Theorem*}{Theorem}
  \newtheorem*{Statement*}{Statement}
  \newtheorem*{Question*}{Question}
  \newtheorem{Proposition}[Definition]{Proposition}
  \newtheorem*{Proposition*}{Proposition}
  
  \newtheorem{Lemma}[Definition]{Lemma}
  \newtheorem{LemmaDefinition}[Definition]{Lemma -- Definition}
  
\theoremstyle{styleRem}
  \newtheorem*{Remark}{Remark}
  \newtheorem*{Remarks}{Remarks}

\newenvironment{myitemize}{\list{--}{\itemsep=0pt}}{\endlist}

\newcommand*\half{\hbox{$\scriptspace\z@\m@th{}^1\mkern-3mu/\mkern-2mu{}_2$}}

\newcommand*\msp{\unskip\penalty\relpenalty\ }

\newcommand*\sseq{}
\let\sseq\subseteq                  
\let\le\leqslant   
\let\ge\geqslant   

\newcommand*\Se[1][n-1]{\mathbb{S}^{#1}}               
\newcommand*\scal[2]{\langle #1\mid\nobreak #2\rangle} 

\newcommand*\card[1]{\left|#1\right|}      
\newcommand*\Card[1]{\bigl|#1\bigr|}       
\newcommand*\Har[1]{\mathcal{H}^{(#1)}}    
\newcommand*\Pol[1]{\mathcal{P}^{(#1)}}    

\newcommand*\abs[1]{\left|#1\right|}   

\newcommand*\NN{\mathbb{N}}   
\newcommand*\ZZ{\mathbb{Z}}   
\newcommand*\RR{\mathbb{R}}   
\newcommand*\CC{\mathbb{C}}   

\newcommand*\Orth{\mathrm{O}}  
\newcommand*\Aut{\mathop{\mathrm{Aut}}} 

\newcommand*\HH{\mathbb{H}}  

\newcommand*\Sh{\mathcal{S}\mkern-3mu\mathit{h}}   

\newcommand*\Th{\mathrm{\theta}_}
\newcommand*\De{\mathrm{\Delta}_}
\newcommand*\Ph{\mathrm{\Phi}}
\newcommand*\Q{\mathrm{Q}}
\newcommand*\R{\mathrm{R}}

\newcommand*\Qp[1]{\mathrm{Q}^{(#1)}}

\newcommand*\refl{\mathrm{r}_}

\newcommand*\reduce{\overline}

\newcommand*\Weyl{\mathrm{W}}

\newcommand*\Z[1]{\mathbf{Z}^{#1}}              
\newcommand*\W[1]{\mathbf{\Gamma}_{\!#1}}       
\newcommand*\varW[1]{\mathbf{\Gamma}'_{\!#1}}   

\newcommand*\A{\mathbf{A}}
\newcommand*\D{\mathbf{D}}
\newcommand*\E{\mathbf{E}}
\newcommand*\F{\mathbf{F}}
\newcommand*\OO{\mathbf{O}}

\newcommand*\Mod{\mathcal{M}} 

\newcommand*\nl{\\[\jot]} 

\newcommand*\setstandardtabs{\kern20pt\=\kern80pt \= \kill}

\makeatletter
\@ifundefined{@@penalty}{\global\let\@@penalty=\penalty}{}
\newcommand*\skippenalty{\global\let\penalty\@@penalty
                \begingroup\afterassignment\endgroup\count@= }
\newcommand*\NOPAGEBREAK{\par\@@penalty\@M\global\let\penalty\skippenalty}

\title{%
Shells of selfdual lattices\\viewed as spherical designs}
\author{Claude \textsc{Pache}\\[\jot]
        \normalsize University of Geneva (Switzerland)\\[\jot]
        \normalsize email: \texttt{Claude.Pache@math.unige.ch}}
\date{27 May 2004}

\begin{document}
\maketitle

\begin{abstract}
  We find out for which $t$
  shells of selfdual lattices and of theirs shadows
  are spherical $t$-designs.
  The method uses theta series of lattices, which are modular forms.
  We analyse fully cubic and Witt lattices,
  as well as all selfdual lattices of rank at most 24.
\end{abstract}


\section*{Introduction}


A nonempty finite subset of a Euclidean sphere provides
approximations for integrals of functions defined on that sphere. 
In this context, such a subset is called a \emph{spherical design} 
and its efficiency is measured by an integer 
$t \ge0$ called its \emph{strength} \cite{DelsarteGoethalsSeidel} 
(precise definitions are given in Section~\ref{SectSphericalDesigns}).
We are interested here in computing (or at least estimating)
the strengths of shells in some \emph{selfdual lattices}.
These problems have natural formulations in terms of
vanishing Fourier coefficients of \emph{modular forms}
which are appropriate \emph{theta series} of the lattices.

The method used in this article was already used for
finding the spherical design strengths of shells of \emph{extremal} (even) lattices.
(Do not confuse with ``extreme lattice''.
The definition of extremal lattices of level~1 is given at the end of Section~\ref{SectClassifModularForms}.) 
See \cite{VenkovExtremal} and \cite[\S16]{VenkovMartinet1} for unimodular case,
and \cite{BachocVenkov} for some other cases.

\medskip

Let $\Lambda$ be a lattice in the standard Euclidian space $\RR^n$;
we denote by $\scal{x}{y}$ the scalar product
of two vectors $x,y \in \RR^n$.
For a positive number $m$, we denote by
\[
  \Lambda_m   := \{ \lambda \in \Lambda \mid \scal{\lambda}{\lambda} = m \}
\]
the \emph{shell} (or \emph{layer}) of norm $m$
(that is to say of radius $\sqrt m$ in the usual sense of Euclidean geometry).

Given a lattice and a positive integer~$t$, we would like to single out
the three following questions:
\begin{enumerate}
  \item[(1)] Is the shell of \emph{minimal norm} a spherical $t$-design?
  \item[(2)] Is \emph{some} shell a spherical $t$-design? 
  \item[(3)] Is \emph{every} shell a spherical $t$-design?
\end{enumerate}

It is quite easy to show that, if any of these question is true for
$t=2$ then the lattice is rational, that is proportional to an integral lattice
(see, e.g., \cite[Chap.~3, \S~1]{MartinetVenkov}).
It is therefore reasonable to restrict the discussion to \emph{integral} lattices.

Question (1) with $t=4$ (or $t=5$ which is equivalent in this case)
asks whether a lattice is \emph{strongly perfect}, 
which is the basic question in \cite{VenkovMartinet1}. 
It is motivated by the classical result of Voronoi
characterising extreme lattices as those which
are ``eutactic and perfect'' \cite{Voronoi}. 
(A lattice is extreme if the density of the corresponding sphere
packing of $\RR^n$ is a local maximum in the space of all lattices.)
Strongly perfect lattices in dimensions $n \le 11$ have been
classified (see \cite{VenkovMartinet1}
for dimensions $n \le 9$ and $n = 11$,
and \cite{NebeVenkov10} for $n = 10$); there are ten
isometry classes, usually denoted by 
$\mathbf{Z}$, $\mathbf{A}_2$, $\mathbf{D}_4$,
$\mathbf{E}_6$, $\mathbf{E}_6^*$,
$\mathbf{E}_7$, $\mathbf{E}_7^*$,
$\mathbf{E}_8$, $\mathit{K}'_{10}$, and ${\mathit{K}'_{10}}^*$
(subscripts indicate dimensions).

Whenever Question (2) has a positive answer,
it is an experimental fact that there exists a
``rather small'' $m$ for which 
$\Lambda_m$ is a spherical $t$-design, 
but we do not know any general result in this direction. 

For the answer to Question~(3) to be positive, 
it is sufficient that the space
$\Har{2j}(\RR^n)^{\Aut(\Lambda)}$
of $\Aut(\Lambda)$-invariant harmonic polynomials 
on $\RR^n$ which are homogeneous of degree $2j$
is reduced to $\{0\}$ for every positive integer $2j\le t$
(see \cite[Thm.~3.12]{GoethalsSeidel81}).

Questions (1) to (3) make sense for various sets associated to a
lattice $\Lambda$, and in particular for
\emph{shadows} of selfdual lattices. 
Recall that, if $\Lambda$ is an odd selfdual lattice with even part
$\Lambda_0 = \{ \lambda \in \Lambda \mid \scal{\lambda}{\lambda} \equiv 0 \bmod2 \}$,
its shadow $\Sh(\Lambda)$ is the complement of $\Lambda$ in
the dual $\Lambda_0^\sharp$ of~$\Lambda_0$.
(Shadows enter naturally the discussion since they
provide efficient tools to compute theta series.)

\medbreak

We denote by $\Z{n}$ the cubic lattice of rank~$n$.
We denote by $\W{n}$ the Witt lattice of rank $n$,
where $n$ is a multiple of~$4$ (see Section~\ref{SectWittLattices} for the definition);
$\W{8}$ is the unique even selfdual lattice of rank~$8$
(also known as the \emph{Korkine-Zolotareff lattice}).
If $R$ is a root system of norm~$2$, we denote by $R^+$ a selfdual lattice of minimal norm~$2$
with $\Lambda_2=R$; it happens that, up to rank~$23$, such a lattice is unique (whenever it exists).
We denote by $k_1\,R_1+\cdots+k_s\,R_s$ the root system whose irreducible components are
$R_1$, \dots, $R_s$ with multiplicity $k_1$, \dots, $k_s$ respectively;
such a root system is called \emph{strongly eutactic}
if all its components have the same Coxeter number.
Recall that the Leech lattice is
the unique even selfdual lattice of rank~$24$ with minimal norm~$4$,
and the shorter Leech lattice is the unique selfdual lattice of rank~$23$
with minimal norm~$3$.
We say that a noncubic selfdual lattice~$\Lambda$ of rank~$n$ has a \emph{long shadow}
if $\sigma(\Lambda)=n-\nobreak8$, where $\sigma(\Lambda)$ is the
minimum norm of a characteristic (or parity) vector of~$\Lambda$.
The definitions of spherical $t$-design and $t\half$-design
are given in Section~\ref{SectSphericalDesigns}.

Theorem~A below is partly a reformulation of known results:
Claims~(i) to (iii) and some cases of Claims (iv)~and~(v)
follow from results on  harmonic polynomials
that are invariant by the action of the automorphism group of the lattice,
as explained above (\cite[Thm.~3.2]{GoethalsSeidel81};
\cite[Examples 7.6~and~7.7]{GoethalsSeidel}).

\begin{Theorem*}[A]\ \NOPAGEBREAK
  \begin{enumerate}
    \item All shells of the Leech lattice
      are spherical $11\half$-designs.
    \item All shells of the shorter Leech lattice and of its shadow
      are spherical $7$-designs. 
    \item All shells of the Korkine-Zolotareff lattice
      are spherical $7\half$-designs.
    \item
      The following special shells of selfdual lattices of rank at most~$24$ and of their shadows
      are spherical $5$-designs:
      \begin{tabbing}\kern40pt \= \kern 100pt \= \kill
        \> $ (\Z{4})_m$                     \> for $m = 2a$,                   \nl
        \> $ (\Z{7})_m$                     \> for $m = 4^a(8b+3)$, $a,b\ge0$, \nl
        \> $ \bigl(\Sh(\Z{16})\bigr)_m$     \> for $m = 4a+2$, $a \ge 1$,      \nl
        \> $ (4\,\A_5)^+_m, (5\,\D_4)^+_m$  \> for $m = 4^a$, $a\ge0$,         \nl
        \> $ (2\,\D_8)^+_m$                 \> for $m = 4a+2$, $a \ge 1$.
      \end{tabbing}
    \item All shells of the following lattices and of their shadows are spherical $3$-designs:
      the cubic lattices $\Z{n}$, the Witt lattices $\W{n}$,
      all selfdual lattices of rank at most~$24$, of minimum~$2$,
      and with strongly eutactic root system
      (see Definition~\ref{DefStronglyEutacticRootSystem};
      that includes all selfdual lattices of minimum~$2$ with long shadow
      and all even selfdual lattices of rank at most~$24$).
      Moreover some other selfdual lattices of rank at most~$24$ have some shells
      which are spherical $3$-designs.
  \end{enumerate}
\end{Theorem*}

Our approach provides a setting for numerical computations;
so we have checked that
no shell of norm at most $1200$ of selfdual lattices up to rank~$24$
or of theirs shadows
is a spherical $t$-design for larger values of $t$ than those indicated in the theorem.

See Theorems \ref{StrengthCubicLattices},
\ref{StrengthWittLattices},
\ref{StrengthEvenUnimodularLattices},
\ref{StrengthLongShadowLattices},
\ref{StrengthLongShadowLatticesBis}, 
and~\ref{StrengthOddUnimodularRank24Lattices}
for details.

\medbreak

Let
\[
  \De{24} = q^2 \prod_{m\ge1}\bigl(1-q^{2m}\bigr)^{24}  =  \sum_{m\ge1} \tau(m)\,q^{2m}
\]
be the generating series of the Ramanujan numbers~$\tau(m)$.
It is a famous conjecture of Lehmer that $\tau(m)$ is 
never zero, and it has been verified for $m\le10^{15}$
\cite[\S~3.3]{Serre85}.
The following Proposition gives a reformulation of that conjecture in terms
of spherical design strengths of shells of the Korkine-Zolotareff lattice
and of the even selfdual lattices of rank~$16$.
It is Proposition~\ref{TheoremOnLehmerConjecture} in our article.

\begin{Proposition*}[B]
  For $m\ge1$, the following are equivalent:\NOPAGEBREAK
  \begin{enumerate}
    \item[(a)] $\tau(m)=0$;
    \item[(b)] the shell of norm $2m$ of the Korkine-Zolotareff lattice
               is an $8$-design (and therefore an $11$-design);
    \item[(c)] the shell of norm $2m$ of any even selfdual lattice~$\Lambda$ of rank~$16$,
               is a $4$-design (and therefore a $7$-design);
  \end{enumerate}
\end{Proposition*}

[Note that, for example, in Condition~(c) above,
the shells of $\Lambda$ are not only $3$-designs, but also $3\half$-designs
(see Definition~\ref{DefAntipodalSphericalHalfDesign});
therefore, if a shell of $\Lambda$ is a $4$-design, then it is a $7$-design.]

Here are other similar equivalences between spherical design strengths of shells of lattices
and vanishing Fourier coefficients of modular forms.
(The definitions of the series $\De8$, $\Th2$, $\Th3$, and $\Th4$
appear in Section~\ref{SectClassifModularForms}.)

\begin{Proposition*}[C]
  Consider any of the following choice of
  a selfdual lattice~$\Lambda$,
  a positive integer~$t$,
  and  a series~$\Theta(z)=\sum_{m\ge1}a_mq^m$ where $q=e^{i\pi z}$.
  \begin{myitemize}
    \item $\Lambda=\Z{n}$ the cubic lattice of rank $n\ge2$,
      $t=3$,
      $\Theta=\De8\,\Th3^n$
      (Section~\ref{SectCubicLattices});
    \item $\Lambda=\W{n}$ the Witt lattice of rank $n\ge12$, $n\equiv0\bmod4$,
      $t=3$,
      $\Theta = \Th2^4\Th3^4\Th4^4\*\bigl(-\Th2^{n-4}+\Th3^{n-4}-\Th4^{n-4}\bigr)$
      (Section~\ref{SectWittLattices});  
    \item $\Lambda$ the Korkine-Zolotareff lattice, of rank~$8$,
      $t=7$,
      $\Theta=\De{24}$
      (Section~\ref{SectEvenLattices24});
    \item $\Lambda$ an even selfdual lattice of rank~$16$,
      $t=3$,
      $\Theta=\De{24}$
      (Section~\ref{SectEvenLattices24});
    \item $\Lambda$ one of the $23$ Niemeier lattice, of rank~$24$,
      $t=3$,
      $\Theta=\Q\,\De{24}$
      (Section~\ref{SectEvenLattices24});
    \item $\Lambda$ the Leech lattice,
      $t=11$,
      $\Theta=\De{24}^2$
      (Section~\ref{SectEvenLattices24});
    \item $\Lambda$ one of the 12 odd selfdual lattices with long shadow and of minimum~2
      (with rank $n=12$ or $14\le n\le22$),
      $t=3$,
      $\Theta=\De8^2\,\Th3^{n-8}$
      (Section~\ref{SectLatticesLongShadow});
    \item $\Lambda$ the shorter Leech lattice, of rank~$23$,
      $t=7$,
      $\Theta=\De8^2\,\Th3^{15}$
      (Section~\ref{SectLatticesLongShadow}).
  \end{myitemize}  
  Then, for each $m\ge1$,
  the shell of $\Lambda$ of norm $m$
  is a spherical $t$-design;
  moreover it is not a spherical $(t+\nobreak1)$-design if and only if $a_m\ne0$.
\end{Proposition*}

It is therefore interesting to look for vanishing coefficients
of the series mentioned above.
For example, we ask the following questions:

\begin{enumerate}
  \item[(1)]
    Consider the series
    \[
      \De{24}^2 = \sum_{m\ge2} a_m q^{2m}, \qquad a_m=\sum_{i+j=m}\tau(i)\,\tau(j),
    \]
    where $\tau(i)$ is the $i$th Ramanujan number.
    Is it true that $a_m\ne0$ for every $m\ge2$~?
  \item[(2)]
    Consider the series
    \[
      \De8\,\Th3^n = \sum_{m\ge1} a_m q^{m}.
    \]
    Is it true that $a_m\ne0$ for every positive integer~$m$ not of the form
    $4^a(8b+3)$, $a,b\ge0$~?
\end{enumerate}

(There are similar questions for the other forms mentioned in Proposition~C.)

We have checked numerically that the answers are positive for $m\le1200$. 


\medbreak

In Sections \ref{SectSphericalDesigns}~to~\ref{SectClassifModularForms},
we recall standard material on spherical designs,
selfdual lattices, theta series, and modular forms.
Section~\ref{SectCalcThetaSeries} gives the form of the theta series
of selfdual lattices;
this is the centre of our analysis.
Some indices of vanishing Fourier coefficients for modular forms
are given in Section~\ref{SectZeroCoeff}.
In Section~\ref{SectRootSystems}, we analyse
the spherical design strengths of root systems of norm~$2$.
Sections \ref{SectCubicLattices}~to~\ref{SectOtherSelfdualLattices}
contain the results
on the strength of shells of some selfdual lattices,
namely the two infinite series of cubic and Witt lattices,
and all selfdual lattices of rank at most~$24$.
Finally, an appendix contains an alternative proof not using modular forms
of the fact
that appropriate shells of $\Z{4}$ and $\Z{7}$ are spherical $5$-designs.

\part{General theory}

\section{Spherical designs}\label{SectSphericalDesigns}

Let $n\ge2$ be an integer, and
let $m$ be a positive real number;
we denote by
\[
  \Se_m:=\{ x\in\RR^n \mid \scal xx = m\}
\]
the sphere of square radius~$m$,
and by $\sigma$ the probability measure on $\Se_m$
invariant under the action of the orthogonal group~$\mathrm{O}(n)$.
A~\emph{spherical design of strength~$t$}, or a \emph{$t$-spherical design},
is a nonempty finite subset $X\sseq\Se_m$ such that
\begin{equation*}
  \frac{1}{\card{X}}\sum_{x\in X} P(x) = \int_{\Se_m} P(y)\,d\sigma(y)
\end{equation*}
for every polynomial form $P$ on~$\RR^n$ of degree at most~$t$
\cite{DelsarteGoethalsSeidel}.

We denote by $\Har{j}(\RR^n)$ the set of homogeneous polynomial forms on~$\RR^n$
of degree~$j$ that are harmonic.
It is classical that $X$ is a spherical $t$-design if and only if
the condition
\begin{equation*}\tag{$C_{j}$}
  \sum_{x\in X} P(x)=0,\qquad \forall P\in\Har{j}(\RR^n)
\end{equation*}
holds for every integer $j$ such that $1\le j\le t$.
It is indeed an immediate consequence of the decomposition
\[
  \Pol{t}(\RR^n) = \Har{t}(\RR^n)\oplus\omega\Pol{t-2}(\RR^n),
\]
where $\Pol{t}(\RR^n)$ is the space of homogeneous polynomial forms on $\RR^n$ of degree~$t$,
and $\omega(x):=\scal{x}{x}$;
see~\cite[\S IX.2]{Vilenkin} and \cite[Chap.~1, \S\S2,~3]{MartinetVenkov}.

In this article, we study spherical designs  which are shells of selfdual lattices
or of their shadows.
These designs are antipodal sets, that is sets~$X$ satisfying $-X=X$;
in this case,
Condition~$(C_j)$ is automatically fulfilled for $j$ odd.
Therefore we use the following reformulation for antipodal sets:

\begin{Definition}\label{DefAntipodalSphericalDesign}
  A nonempty finite antipodal subset
  $X\sseq\Se_m$ is a \emph{spherical $(2s+\nobreak1)$-design}
  (or, equivalently, a spherical $2s$-design) if the condition
  \begin{equation*}\tag{$C_{2j}$}%
    \sum_{x\in X} P(x)=0,\qquad \forall P\in\Har{2j}(\RR^n)
  \end{equation*}
  holds for every even integer~$2j$ such that $2\le 2j\le 2s$.
\end{Definition}

Some antipodal spherical $(2s+1)$-designs, although not satisfying Condition $(C_{2s+2})$,
do verify Condition $(C_{2s+4})$. Therefore, we define \cite{VenkovExtremal}: 

\begin{Definition}\label{DefAntipodalSphericalHalfDesign}
  A nonempty finite antipodal subset  $X\sseq\Se_m$ is a
  \emph{spherical $(2s+\nobreak1+\nobreak\frac{1}{2})$-design}
  if it verifies condition $(C_{2j})$
  for $2\le 2j\le 2s$ and $2j=2s+4$.
\end{Definition}

\section{Selfdual lattices and shadows}\label{SectLattices}

For a subset $A\sseq\RR^n$ and a positive real number $m>0$,
the \emph{shell} (or \emph{layer})  of norm~$m$ of~$A$ is
\[
  A_m:=\{ x\in\Lambda \mid \scal xx =m\} = A \cap \Se_m.
\]
A \emph{lattice} of rank~$n$ is a discrete $\ZZ$-submodule of $\RR^n$
which spans $\RR^n$ as $\RR$-module.
Two lattices are \emph{equivalent} if there exists an orthogonal linear transformation
which sends one lattice onto the other;
we often consider two equivalent lattices as being the same lattice.
We define the \emph{minimum} (or the \emph{minimal norm}) of $\Lambda$ as
\[
  \min(\Lambda) := \min \{m>0\mid\Lambda_m\ne\emptyset\}.
\]
The \emph{dual} of a lattice $\Lambda$ is the lattice
\[
  \Lambda^\sharp:=\{y\in\RR^n \mid \scal{y}{x}\in\ZZ\ \forall x\in\Lambda\}.
\]
The lattice $\Lambda$ is called \emph{integral} if $\Lambda\sseq\Lambda^\sharp$, that is if
$\scal xy \in\ZZ$ for all $x,y\in\Lambda$.
An integral lattice~$\Lambda$ is called \emph{even} if
$\scal{x}{x}\in 2\ZZ$ for all $x\in\Lambda$;
it is called \emph{odd} otherwise.
An integral lattice is called \emph{selfdual} (or \emph{unimodular})
if $\Lambda^\sharp=\Lambda$.

For $A\sseq\RR^a$ and $B\sseq\RR^b$ we set
$A\oplus B:= \{(x,y)\in\RR^{a+b} \mid x\in A,\msp y\in B\}$.
It is easily checked that any integral lattice $\Lambda$ of rank~$n$
is of the form
\[
  \Lambda\simeq\Z{p}\oplus L,
\]
where $L$ is an integral lattice of rank $N=n-p$
and of minimum at least~$2$,
and where
\[
   \Z{p}:=\{(x_1,\dots,x_p) \mid x_i\in\ZZ\}\sseq\RR^p
\]
is the cubic lattice of rank~$p$.
Note that $L$ is selfdual if and only if $\Lambda$ is selfdual.

The \emph{shadow} $\Sh(\Lambda)$ of a selfdual lattice $\Lambda$ is defined as follows:
If $\Lambda$ is even, we set $\Sh(\Lambda):=\Lambda$. Otherwise let
\begin{equation*}
  \Lambda_0:=\{x\in\Lambda \mid \scal xx \equiv 0\bmod2\},
\end{equation*}
which is an even sublattice of $\Lambda$ of index~$2$;
therefore $\Lambda$ is a sublattice of $\Lambda_0^\sharp$ of index~$2$.
We set
\[
  \Sh(\Lambda):=\Lambda_0^\sharp \setminus \Lambda.
\]
An alternative description of the shadow is
$\Sh(\Lambda)=\{x/2\mid x$ is a characteristic vector of~$\Lambda\}$,
where a \emph{characteristic vector} (or \emph{parity vector}) of $\Lambda$
is a vector~$x\in\Lambda$ such that $\scal{x}{y}\equiv\scal{y}{y}\penalty\binoppenalty\bmod2$
for all $y\in\Lambda$.

For a selfdual lattice~$\Lambda$, we define:
\[
  \sigma(\Lambda) := 4\min \{\scal xx \mid x\in\Sh(\Lambda)\},
\]
which is a nonnegative integer.
(It is the minimal norm of the characteristic vectors
of $\Lambda$.)
We have $\sigma(\Lambda)=0$ if and only if $\Lambda$ is even.

It is easily checked that, if $\Lambda'$ and $\Lambda''$ are
selfdual lattices, then $\Lambda'\oplus\Lambda''$ is a selfdual lattice, and
\begin{gather*}
  \Sh(\Lambda'\oplus\Lambda'') = \Sh(\Lambda')\oplus\Sh(\Lambda''), \\
  \sigma(\Lambda'\oplus\Lambda'') = \sigma(\Lambda')+\sigma(\Lambda'').
\end{gather*}

The following facts are well-known;
a proof using modular forms appears at the end of Section~\ref{SectClassifModularForms}:

\begin{Proposition}\label{BasicFactsShadow}
  Let $\Lambda\sseq\RR^n$ be a selfdual lattice. Then
  \begin{enumerate}
    \item for every $x\in\Sh(\Lambda)$, we have $4\scal xx \equiv n \bmod 8$;
    \item there exists a nonnegative integer $k$ such shat $\sigma(\Lambda)=n-8k$.
  In particular, if $\Lambda$ is even,
  then $n\equiv 0\bmod 8$;
    \item we have $\sigma(\Lambda)=n$ if and only if $\Lambda\simeq\Z{n}$,
  \end{enumerate}
\end{Proposition}

The characterisation of $\Z{n}$ given in Claim~(iii)
is due to Elkies \cite{Elkies0}.

Note that, in the decomposition $\Lambda\simeq\Z{p}\oplus L$,
if $\sigma(\Lambda)=n-8k$, then $\sigma(L)=N-8k$,
where $n$ is the rank of $\Lambda$,
and $N$ is the rank of~$L$.

The list of selfdual lattices of rank at most~$24$ can be found in
\cite[Chap.\ 16~and~17]{ConwaySloane}
and \cite{Bacher}.

\section{Theta series}\label{SectThetaSeries}

Let us now introduce the tools for analysing the
spherical design strengths of shells of selfdual lattices.
Let
\[
    \HH:=\{z\in\CC \mid \Im z>0\}
\]
be the Poincar\'e half-plane.
Recall that a holomorphic function $f:\HH\to\CC$ is
\emph{bounded at infinity} if there exists $r>0$ such that
$\{f(z)\mid \Im z\ge r\}$ is bounded.

\begin{LemmaDefinition}\label{LemThetaSeries}
  Let $A$ be a nonempty subset of~$\RR^n$
  for which there exists $\delta>0$ such that
  $\abs{x-y}\ge\delta$ for all distinct $x,y\in A$.
  Let~$P$ be a polynomial form on $\RR^n$.
  Then the series
  \[
    \Theta_{A,P}(z) := \sum_{x\in A} P(x) \, e^{i\pi z\scal xx},
     \qquad z\in\HH,
  \]
  converges absolutely to
  a function on $\HH$ that is holomorphic and bounded at infinity.
  It is the \emph{theta series} of $A$ weighted by~$P$.

  Let us assume moreover that there exists a real number $\alpha>0$
  such that $\alpha\scal xx \in 2\ZZ$ for every $x\in A$.
  Then we have
  \[
    \Theta_{A,P}(z+\alpha)=\Theta_{A,P}(z),\qquad \forall z\in\HH.
  \]
\end{LemmaDefinition}

\begin{proof}
  The condition on the distance of two distinct points of~$A$
  implies that there exist a constant~$C>0$ such that
  for every $r\ge 0$ the set
  \[
    \{ x\in A \mid r \le \scal xx \le r+1 \}
  \]
  contains at most $C\,r^{n-1}$ points;
  the absolute convergence of $\Theta_{A,P}$ follows.
  The second claim of the lemma is straightforward.
\end{proof}

\begin{Remarks} \
  \NOPAGEBREAK
  \begin{enumerate}
    \item
      For a holomorphic function $F:\HH\to\CC$ verifying
      $F(z+\alpha)=F(z)$ for all $z\in\HH$ and for some $\alpha>0$,
      the condition to be bounded at infinity
      is equivalent to the condition that there exists a Fourier expansion
      of the form
      \[
         F(z) = \sum_{\!m\in 2\alpha^{-1}\NN\!} a_m e^{i\pi zm}
      \]
      which converges for $\Im z$ sufficiently large.
      (We use $\NN=\{0,1,2,\ldots\}$.)
      In this case, $F$ is said to be \emph{holomorphic at infinity}.
    \item
      For a real number~$m$ and if $z\in\HH$ is understood, we write
      \[
        q^m:=e^{i\pi z m},\quad \text{where $z\in\HH$.}
      \]
      Thus the theta series of~$A$ weigted by~$P$ can be written
      \[
        \Theta_{A,P}(z) = \sum_{x\in A}P(x)\,q^{\scal xx} = \sum_{\!m\in 2\alpha^{-1}\NN\!} a_m^{(P)} q^m,
        \qquad\text{where } a_m^{(P)}:=\sum_{x\in A_m} P(x).
      \]
    \item
      The classical theta series of~$A$ is
      \[
        \Theta_A:=\Theta_{A,1} = \sum_{m\ge0} \card{A_m}q^m.
      \]
    \item Let $A\sseq\RR^a$ and $B\sseq\RR^b$,
      and let $P$, respectively $Q$, be polynomial forms on $\RR^a$, respectively $\RR^b$. Then
      \[
        \Theta_{A\oplus B,\, PQ} = \Theta_{A,P}\,\Theta_{B,Q}.
      \]
  \end{enumerate}
\end{Remarks}
  
Now we reformulate the condition of being a spherical design using theta series.

\begin{Lemma}\label{CharactShellIsDesign}
  Let $A$ be a nonempty subset of~$\RR^n$
  such that there exist $\delta>0$ verifying
  $\abs{x-y}\ge\delta$ for all distinct $x,y\in A$.
  Then, for~$m>0$,
  the shell $A_m$ is a spherical $t$-design or is empty
  if and only if
  \[
    \text{$a_m^{(P)}=0$ for every $P\in\Har{2j}(\RR^n)$, $1\le 2j\le t$,}
  \]
  where $a_m^{(P)}$ are the Fourier coefficients of the theta series
  \[
    \Theta_{A,P}(z) = \sum_{m} a_m^{(P)} q^m.
  \]
\end{Lemma}
\begin{proof}
  This follows directly from the definitions.
  Note that the Fourier coefficients $a_m^{(P)}$ are retrieved from $\Theta_{A,P}$
  by the formula
  $$
    a_m^{(P)}=\lim_{R\to\infty} \frac{1}{2R} \int_{-R}^{R} e^{-i\pi m(x+iy_0)}\Theta_{A,P}(x+iy_0)\,dx,
  $$
  valid for any $y_0>0$.
\end{proof}

\section{Modular forms}\label{SectModularForms}

The interest of theta series of lattices is that they have the following property,
which will help to compute them, at least if the lattice is selfdual.
This proposition is a direct consequence of the Poisson Summation Formula;
see for example \cite[Prop.~3.1, p.~87]{Ebeling}.

\begin{Proposition}\label{CorFormulaPoisson}
  Let $\Lambda\sseq\RR^n$ be a lattice, and let $P:\RR^n\to\CC$
  be a harmonic polynomial of degree~$2j$.
  Then
  \begin{equation*}
    \Theta_{\Lambda^\sharp,P}(z) = (\det\Lambda)^{1/2} (-1)^j (i/z)^{n/2+2j} \Theta_{\Lambda,P}(-1/z) .
  \end{equation*}
  (For the power of $i/z$, we use the principal branch;
  observe that $-\pi/2\le\arg(i/z)\le\pi/2$ for $z\in\HH$.)
\end{Proposition}

In the case of selfdual lattices, the latter formula gives a relation
between $\Theta_{\Lambda,P}(z)$ and $\Theta_{\Lambda,P}(-1/z)$.
To be more precise, let us give the following definitions.
Recall that $\NN=\{0,1,2,\ldots\}$.

\begin{Definition}\label{DefModularForm}\ \NOPAGEBREAK
  \begin{enumerate}
    \item
      Let $\lambda\in\{1,2\}$, $\omega\in\frac12\NN$, and $\epsilon\in\{+,-\}$.
      A \emph{modular form of signature $(\lambda,\omega,\epsilon)$}
      is a holomorphic and holomorphic at infinity function $f:\HH\to\CC$
      that verifies
      \begin{gather*}
        f(z+\lambda) = f(z),  \\
        f(-1/z) = \epsilon (z/i)^\omega f(z),
      \end{gather*}
      for all $z\in\HH$.
      The number $\omega$ is the \emph{weight} of the form.
      We denote by $\Mod_\omega^{\lambda,\epsilon}$ the vector space of
      modular forms of signature $(\lambda,\omega,\epsilon)$;
      note that $\Mod_\omega^{1,\epsilon}\sseq\Mod_\omega^{2,\epsilon}$.
      We set:
      \[
        \Mod^{\lambda,\epsilon} := \bigoplus_{\!\omega\in(1/2)\NN\!} \Mod_\omega^{\lambda,\epsilon},
        \qquad
        \Mod^{\lambda} := \Mod^{\lambda,+} \oplus \Mod^{\lambda,-}.
      \]
      Observe that $\Mod^{\lambda}$ and $\Mod^{\lambda,+}$ are algebras, graded by the weight.

    \item
        A modular form~$f$ is \emph{parabolic} if
        \[
          \lim_{\Im z\to +\infty} f(z) = 0.
        \]
  \end{enumerate}
\end{Definition}

Nonzero modular forms of signature $(1,\omega,\epsilon)$ exist only for
weights $\omega\in 4\NN$ if $\epsilon=\mathord{+}$,
and for weights $\omega\in 4\NN+2$ if $\epsilon=\mathord{-}$.
Indeed, let $S:z\mapsto -1/z$ and $T:z\mapsto z+1$. We have $STSTST=\mathrm{id}_{\HH}$,
therefore $f(STSTSTz)=f(z)$ for every $z\in\HH$.      
On the other hand, it is straightforward to check that, for $f\in\Mod_\omega^{1,\epsilon}$,
\[
  f(STSTSTz) = \epsilon\, i^\omega f(z), 
\]
where $i$ must be taken with argument $\pi/2$. Thus $\epsilon\, e^{i\omega\pi/2}=1$ if $f\ne0$.

Nonzero modular forms of signature $(2,\omega,\epsilon)$ exist for every
weight $\omega\in\frac12\NN$;
for example, according to Proposition~\ref{SeriesThetaSelfdualLattice} below,
$\Th3^n=\Theta_{\Z{n}}$ is a modular form of parameters $(1,n/2,+)$.

Modular forms of signature $(1,\omega,\epsilon)$,
are the classical modular forms for $\mathrm{SL}(2,\ZZ)$.
Modular forms of signature
$(2,\omega,\epsilon)$
are modular forms for the
subgroup of index~3 in $\mathrm{SL}(2,\ZZ)$ whose elements are the matrices that reduce
to $\left(\begin{smallmatrix}1&0\\0&1\end{smallmatrix}\right)$
or $\left(\begin{smallmatrix}0&1\\1&0\end{smallmatrix}\right)$
modulo~$2$ (sometimes noted $G(2)$ or $\Gamma_V(2)$).

Proposition~\ref{CorFormulaPoisson} implies immediately:

\begin{Proposition}\label{SeriesThetaSelfdualLattice}
  Let $\Lambda\sseq\RR^n$ be a selfdual lattice, and let $P\in\Har{2j}(\RR^n)$.
  Then
  \begin{align*}
    \Theta_{\Lambda,P} \in \Mod_{n/2+2j}^{2,+} &\quad\text{if $2j\equiv 0\bmod 4$}, \\
    \Theta_{\Lambda,P} \in \Mod_{n/2+2j}^{2,-} &\quad\text{if $2j\equiv 2\bmod 4$}.
  \end{align*}
  If moreover $\Lambda$ is even, then 
  \begin{align*}
    \Theta_{\Lambda,P} \in \Mod_{n/2+2j}^{1,+} &\quad\text{if $2j\equiv 0\bmod 4$}, \\
    \Theta_{\Lambda,P} \in \Mod_{n/2+2j}^{1,-} &\quad\text{if $2j\equiv 2\bmod 4$}.
  \end{align*}
  Moreover, if $2j>0$ then $\Theta_{\Lambda,P}$ is parabolic.
\end{Proposition}

\medbreak
Let us now look at theta series of shadows of selfdual lattices.

\begin{Definition}\label{DefShadow}
  Let $f\in\Mod_\omega^{2,\epsilon}$.
  The \emph{shadow} of $f$ is the function $\Sh f:\HH\to\CC$
  defined by
  \[
    \Sh f (z) := (i/z)^\omega f(-1/z+1).
  \]
\end{Definition}

\begin{Proposition}\label{ThetaSeriesOfShadow}
  Let $\Lambda$ be a selfdual lattice, and let $P\in\Har{2j}(\RR^n)$.
  Then
  \[
    \Theta_{\Sh(\Lambda),P} = (-1)^j \Sh\Theta_{\Lambda,P}.
  \]
\end{Proposition}

\begin{proof}
  If $\Lambda$ is even, this follows immediatly from $\Sh(\Lambda)=\Lambda$
  and from Proposition~\ref{SeriesThetaSelfdualLattice}.
  If $\Lambda$ is odd, we observe that
  \[
    2\Theta_{\Lambda_0,P}(z)=\Theta_{\Lambda,P}(z)+\Theta_{\Lambda,P}(z+1),
  \]
  where $\Lambda_0$ is the even sublattice of $\Lambda$, of index~$2$.
  Then we apply Proposition~\ref{CorFormulaPoisson} to $\Lambda_0$
  and $\Lambda$ to obtain, for $P\in\Har{2j}(\RR^n)$,
  \begin{align*}
    \Theta_{\Sh(\Lambda),P}
      &= \Theta_{\Lambda_0^\sharp,P}(z) - \Theta_{\Lambda,P}(z)  \\
      &= 2(-1)^j(i/z)^{n/2+2j} \Theta_{\Lambda_0,P}(-1/z)
         - (-1)^j(i/z)^{n/2+2j}\Theta_{\Lambda,P}(-1/z)  \\
      &= (-1)^j\Sh\Theta_{\Lambda,P}(z).
    &\rlap{\kern-.5em\null\qedhere}
  \end{align*}
\end{proof}

\section{A theorem of classification for modular forms}\label{SectClassifModularForms}

In order to use Proposition~\ref{SeriesThetaSelfdualLattice},
we need a description of the vector spaces $\Mod^{\lambda,\epsilon}_\omega$
defined in the previous section.
There are very nice results of classification that we recall in this section.

Before giving examples of modular forms, it is convenient
(although not essential) to formulate a definition
of the weight which applies to a larger class of functions than modular forms:

\begin{Definition}\label{DefWeight}
  Let $\omega\in\RR$.
  Let $f$ and $g$ be two meromorphic functions on~$\HH$.
  Assume that there exist $\alpha$ and $\beta>0$ such that
  \begin{gather*}
    f(z+\alpha) = f(z), \\
    g(z+\beta) = g(z),  \\
    f(-1/z) = (z/i)^\omega g(z),
  \end{gather*}
  for all $z\in\HH$.
  (Note that the last condition is symmetric in $f$ and $g$.)
  We say that $f$ and $g$ are of \emph{weight}~$\omega$.
\end{Definition}

When $f$ is a modular form of signature $(\lambda,\omega,\epsilon)$
(Definition~\ref{DefModularForm}),
the assumptions of the Definition above
are verified for $\alpha=\beta=\lambda$ and $g=\epsilon\,f$,
and so $f$ is of weight~$\omega$ as expected.
However, some functions on~$\HH$ that are holomorphic and holomorphic at infinity,
although not being modular forms, do have a weight in the sense of Definition~\ref{DefWeight}.
For example, we have $\Th4(-1/z)=(z/i)^{1/2}\Th2(z)$, where $\Th4$ and $\Th2$
are the two periodic functions in~$z$ defined below.
Therefore $\Th2$ and $\Th4$ are both of weight~$1/2$; but they are not modular forms.

It is easily checked that the algebra of meromorphic functions~$f$ on~$\HH$
satisfying the assumptions of Definition~\ref{DefWeight} with
some meromorphic function~$g$ and some \emph{rational} positive numbers $\alpha$ and~$\beta$
is an algebra graded by the weight.

We give here a list of functions
that have a weight in the sense of the Definition~\ref{DefWeight}.
Recall that we write $q^m:=e^{i\pi zm}$.

\begingroup
  \renewcommand*\nl{\displaybreak[0]\\[\jot]}
  \begin{align*}
    \Th2(z)    &= \sum_{\!\!\!\!\!m\in\ZZ+1/2\!\!\!\!\!} q^{m^2} = 2q^{1/4}\bigl(1+q^2+O(q^6)\bigr)
                                                                         && \text{of weight $1/2$,} \nl
    \Th3(z)    &= \sum_{m\in\ZZ} q^{m^2}     = 1+2q+2q^4+O(q^9)          && \text{of weight $1/2$,} \nl
    \Th4(z)    &= \sum_{m\in\ZZ} (-q)^{m^2}  = 1-2q+2q^4+O(q^9)          && \text{of weight $1/2$,} \nl
     & \qquad\Th2^4+\Th4^4=\Th3^4 \nl
    \Ph(z)     &= \Th4^4(z)-\Th2^4(z)        = 1-24q+24q^2-96q^3+O(q^4) && \text{of weight $2$,} \nl
    \De8(z)    &= \frac{1}{16}\Th2^4(z)\,\Th4^4(z)
                                             = q-8q^2+28q^3+O(q^4)       && \text{of weight $4$,} \nl
    \Q(z)      &= \Th3^8(z)-16\De8(z)        = 1+240q^2+2160q^4+O(q^6)   && \text{of weight $4$,} \nl
    \R(z)      &= \Ph(z)\bigl(\Th3^8(z)+8\De8(z)\bigr)
                                             = 1-504q^2+O(q^4)           && \text{of weight $6$,} \nl
    \De{24}(z) &= \Th3^8(z)\, \De8^2(z)       = q^2-24q^4+252q^6+O(q^8)  && \text{of weight $12$.}
  \end{align*}
\endgroup
The coefficients of the latter form
\[
  \De{24}(z) = \sum_{m\ge1} \tau(m) q^{2m}
\]
are the celebrated Ramanujan numbers.
All the functions listed above can be expressed in terms of the Jacobi theta function
\[
  \Th3(\xi\mid z) = \sum_{m\in\ZZ} e^{2i\xi m + i\pi z m^2}, \quad z\in\HH;
\]
indeed, we have
\[
  \Th2(z) = e^{i\pi z/4}\Th3(\pi z/2 \mid z), \qquad
  \Th3(z) = \Th3(0 \mid z), \qquad
  \Th4(z) = \Th3(\pi/2 \mid z).
\]

Let $\Mod$ be the algebra generated by $\Th2$, $\Th3$ and $\Th4$.
It is an algebra graded by the weight,
where the weight ranges over the set of nonnegative half-integer.
The elements of weight zero are the constants.

The definition of the shadow of a modular form (Definition~\ref{DefShadow})
carries over
to functions having a weight in the sense of Definition~\ref{DefWeight},
and provides an endomorphism of the graded algebra~$\Mod$.
We give here shadows of some functions in~$\Mod$.
\begingroup
  \renewcommand*\nl{\displaybreak[0]\\[\jot]}
  \begin{align*}
    \Sh\Th2 = \sqrt{i} \Th4, \qquad
    \Sh\Th4 = \Th3, \qquad
    \Sh\Th3 = \Th2 , \nl
    \Sh\De8(z) = -\frac{1}{16}+q^2-7q^4+O(q^6),\nl
    \Sh\Ph(z) = 2+48q^2+48q^4+O(q^6),\nl
    \Sh\Q=\Q, \qquad \Sh\R=-\R, \qquad \Sh\De{24}=\De{24}.
  \end{align*}
\endgroup
These formulae (and many other useful identities)
can be found in \cite[Chap.~4, \S4.1]{ConwaySloane}.

All modular forms of signature $(1,\omega,\pm)$ or $(2,\omega,\pm)$
are elements of~$\Mod$.
More precisely:

\begin{Theorem}\label{ClassModularForms}
  \begin{align*}
     &\Mod^{2,+} =      \CC[\Th3,\De8],     &&\Mod^{1,+} =      \CC[\Q,\De{24}], \\
     &\Mod^{2,-} = \Ph\,\CC[\Th3,\De8],     &&\Mod^{1,-} = \R \,\CC[\Q,\De{24}].
  \end{align*}
\end{Theorem}
For a proof of this deep result, see \cite{Rankin},
Section~6.1 for $\Mod^{1,\pm}$, and Section~7.1 for $\Mod^{2,\pm}$.
The notation of \cite{Rankin} are:\\
\indent
  $\{\Gamma(1),\omega\}$ for $\Mod_\omega^{1,+}$ when $\omega\in4\NN$,
   and for $\Mod_\omega^{1,-}$ when $\omega\in4\NN+2$;\\
\indent
  $\{\Gamma_V(2),\omega,v_{2\omega}\}$
  for $\Mod_\omega^{2,+}$, $\omega\in\frac12\NN$,
  and\\
\indent
  $\{\Gamma_V(2),\omega,v^*_{2\omega}\}$
  for $\Mod_\omega^{2,-}$, $\omega\in\frac12\NN$.\\
Moreover modular forms are called there \emph{integral modular forms}.

\medbreak
Here is a first application of this theorem:

\begin{proof}[Proof of Proposition~\ref{BasicFactsShadow}]
  For brevity, we write $\Theta_{\Lambda}$, $\De8$, etc.\
  instead of  $\Theta_{\Lambda}(z)$, $\De8(z)$, etc.
  By Theorem~\ref{ClassModularForms} and Proposition~\ref{SeriesThetaSelfdualLattice},
  for $\Lambda$ a selfdual lattice,
  the theta series
  \[
    \Theta_{\Lambda} = 1 + \sum_{m\ge1}\card{\Lambda_m}q^m
  \]
  is a polynomial in $\Th3$ and $\De8$ and is of weight~$n/2$;
  therefore it is of the form
  \[
    \Theta_{\Lambda}=\Th3^n + c_1\Th3^{n-8}\De8 + \cdots + c_k\Th3^{n-8k}\De8^k
      = 1 + \card{\Lambda_1} q + \card{\Lambda_2} q^2+\cdots ,
  \]
  with $c_k\ne0$.
  (There is no coefficient in front of $\Th3^n$,
  because the series expansion in~$q$ begins with $1+\cdots$.)
  By Proposition~\ref{ThetaSeriesOfShadow},
  the theta series of  the shadow is then
  \begin{align*}
    \Theta_{\Sh(\Lambda)}
      &= \sum_{m}\card{\Sh(\Lambda)_m}q^m
       = \Sh\Theta_\Lambda \\
      &= \Th2^n + c_1\Th2^{n-8}\Sh\De8 + \cdots + c_k\Th2^{n-8k}\Sh\De8^k \\
      &=  d_1 q^{(n-8k)/4} +  d_2 q^{(n-8k+8)/4}+\cdots,
  \end{align*}
  where $d_1=(-1/16)^k 2^{n-8k}c_k\ne0$. So, $\Sh(\Lambda)_m\ne\emptyset$ implies
  $4m\equiv n \allowbreak\bmod8$, and we have $\sigma(\Lambda)=n-8k$.
  In particular, if $\Lambda$ is even, we have $\sigma(\Lambda)=0=n-n$,
  therefore $n$ is a multiple of~$8$.

  Finally, it follows from the decomposition $\Lambda\simeq\Z{p}\oplus\ L$,
  with $L$ a selfdual lattice of minimum at least~$2$,
  that 
  $\card{\Lambda_1}= 2p\le 2n$
  with equality if and only if $\Lambda\simeq\Z{n}$.
  On the other hand, if $\sigma(\Lambda)=n$ then
  $\Theta_\Lambda=\Th3^n = 1+2nq+O(q^2)$, and therefore $\card{\Lambda_1}=2n$.
  Thus $\sigma(\Lambda)=n$ implies $\Lambda\simeq\Z{n}$.
\end{proof}

\begin{Remark}
  Let $\omega>0$ be a multiple of~$4$.
  From Theorem~\ref{ClassModularForms},
  we have $\dim\Mod^{2,+}_\omega=m+\nobreak1$, where $m=[\omega/12]$,
  and there exists a unique $\Theta_\omega\in\Mod^{2,+}_\omega$
  such that
  \[
    \Theta_\omega(q) = 1 + a_{2m+2} q^{2m+2} +O(q^{2m+4}).
  \]
  Moreover, it is known that $a_{2m+2}>0$.
  Such a theta series is called \emph{extremal},
  and an even selfdual lattice $\Lambda$ of rank $2\omega$
  with $\Theta_\Lambda=\Theta_\omega$ is called \emph{extremal}
  (or, more precisely, \emph{extremal of level~1}).
  
  The method used in our paper also apply to extremal lattices:
  see \cite{VenkovExtremal} and \cite[\S16]{VenkovMartinet1}.
  See also \cite{BachocVenkov}, where the case
  of non-selfdual extremal lattices is also treated. 
  
  See \cite{ScharlauSchulze-Pillot} for more informations on extremal lattices.
\end{Remark}

\section{Computing the theta series of a selfdual lattice}\label{SectCalcThetaSeries}

  We can now give more precisely the general form of the theta series of a selfdual lattice.

  \begin{Proposition}\label{CalcThetaSeries}
    Let $\Lambda\sseq\RR^n$ be a selfdual lattice with $\sigma(\Lambda)=n-8k$.
    Let us write 
    \[
      \Lambda=\Z{p}\oplus L, \qquad L\sseq\RR^{N},
    \]
    where $L$ is of minimum at least~$2$.
    Then there exist $c_i\in\ZZ$ such that
    \begin{gather*}
      \Theta_\Lambda = \Th3^n + c_1\De8\Th3^{n-8} + c_2 \De8^2\Th3^{n-16} + \cdots+ c_k\De8^k\Th3^{n-8k}.
    \end{gather*}
    Some values of $c_i$ are:
    \begin{align*}
      c_1 &= -2N \\
      c_2 &= (h - 46+2N)N,\quad\text{where $h:=\card{L_2}/N$} \\
      c_k &= (-1)^k 2^{-n+12k} \card{\Sh(\Lambda)_{(n-8k)/4}}.
    \end{align*}
  \end{Proposition}

  \begin{proof}
    As in the proof of Proposition~\ref{BasicFactsShadow} given at the end of the previous Section,
    we can write
    \[
      \Theta_{\Lambda}= \Th3^n + c_1\Th3^{n-8}\De8 + \cdots + c_k\Th3^{n-8k}\De8^k,
    \]
    with $c_i\in\RR$ and $c_k\ne0$, and $k$ is given by $\sigma(\Lambda)=n-8k$.
    The first coefficients of the Fourier expansion in $q=e^{i\pi z}$ are
    \begin{align*}
        \Theta_{\Lambda} &= (1+2q+2q^4+\cdots)^n + c_1(1+2q+\cdots)^{n-8}(q-8q^2+\cdots) \\
                         & \qquad + c_2(1+\cdots)^{n-16}(q+\cdots)^2+\cdots \\
                         &= \bigl(1+2nq+2n(n-1)q^2+\cdots\bigr) + c_1\bigl(q+(2n-24)q^2+\cdots\bigr) + c_2q^2 + \cdots \\
                         &= 1+ (2n + c_1)q + \bigl(2n(n-1)+(2n-24)c_1+c_2\bigr) q^2 + \cdots.
    \end{align*}
    In particular, we have $\card{\Lambda_1} = 2n + c_1$. Since $\card{\Lambda_1}=2p=2(n-N)$,
    we have $c_1=-2N$. The second coefficient is then
    \begin{align*}
      \card{\Lambda_2} &= 2n(n-1)+(2n-24)c_1+c_2 \\
                       &= 2(p+N)(p+N-1) - (2p+2N-24)2N + c_2 \\
                       &= 2p(p-1) +(46-2N)N + c_2.
    \end{align*}
    Since $\card{\Lambda_2}=2p(p-1)+\card{L_2}$, this gives $c_2= (h-46+2N)N$ where $h=\card{L_2}/N$.

    By induction on~$i$, it is clear that $c_i$ is integral.

    Now, the theta series of the shadow is
    \begin{align*}
        \Theta_{\Sh(\Lambda)}=  \Sh\Theta_\Lambda
          &= c_k\Th2^{n-8k}\Sh\De8^k + \cdots + c_1\Th2^{n-8}\Sh\De8 + \Th2^n \\
          &= c_k\bigl(2q^{1/2} + \cdots\bigr)^{n-8k} \bigl(-1/16+\cdots\bigr)^{k} + \cdots \\
          &= c_k 2^{n-8k} (-1)^k 2^{-4k} q^{(n-8k)/4} + \cdots \\
          &= (-1)^k 2^{n-12k} c_k q^{(n-8k)/4} + \cdots;
    \end{align*}
    this shows that $\card{\Sh(\Lambda)_{(n-8k)/4}}=(-1)^k 2^{n-12k} c_k$.
  \end{proof}

  \begin{Proposition}\label{CalcThetaSeriesBis}
    Let $\Lambda\sseq\RR^n$ be a selfdual lattice with $\sigma(\Lambda)=n-8k$
    and of minimum~$m$.
    \begin{enumerate}
      \item For every even positive integer~$j$,
          there exist linear forms $c_i:\Har{2j}(\RR^n)\to\CC$ such that
          \[
            \Theta_{\Lambda,P} = \sum_{i=m}^{k+j/2} c_i(P)\,\De8^i\Th3^{n+4j-8i},
            \qquad\forall P\in\Har{2j}(\RR^n).
          \]
          In particular, if $m>k+j/2$,
          then $\Theta_{\Lambda,P}=0$ for every $P\in\Har{2j}(\RR^n)$.
      \item For every odd positive integer~$j$,
          there exist linear forms $c_i:\Har{2j}(\RR^n)\to\CC$ such that
          \[
            \Theta_{\Lambda,P} = \sum_{i=m}^{\!\!k+(j-1)/2\!\!} c_i(P)\,\Ph\,\De8^i\Th3^{n+4j-2-8i},
            \qquad\forall P\in\Har{2j}(\RR^n).
          \]
          In particular, if $m>k+(j-1)/2$,
          then $\Theta_{\Lambda,P}=0$ for every $P\in\Har{2j}(\RR^n)$.
    \end{enumerate}
  \end{Proposition}

  \begin{proof}
    We prove only Claim~(i), since the proof of Claim~(ii) is similar.
    
    Let $j$ be an even positive integer.
    By Proposition~\ref{SeriesThetaSelfdualLattice},
    for each $P\in\Har{2j}(\RR^n)$,
    $\Theta_{\Lambda,P}$ is a parabolic form of weight~$n/2+2j$.
    By Theorem~\ref{ClassModularForms},
    it is of the form
    \[
       \Theta_{\Lambda,P} = \sum_{i\ge 1} c_i(P)\,\De8^i\Th3^{n+4j-8i}.
    \]
    Since $\Theta_{\Lambda,P}$ is linear in~$P$,
    the coefficients $c_i$ are also linear in~$P$.

    Let us now suppose that $c_i\not\equiv0$ for some index~$i$,
    and let $a$ [respectively~$b$] be the smallest [respectively the largest] index~$i$
    such that $c_i\not\equiv 0$.
    It remains to prove that
    $a\ge m$ and $b\le k+j/2$.
    We have
    \begin{align*}
       \Theta_{\Lambda,P} &= c_a(P)\,\Th3^{n+4j-8a}\De8^a+\cdots \\
                          &= c_a(P)\,q^a + \cdots,       
    \end{align*}
    hence $c_a(P)=\sum_{x\in\Lambda_a}P(x)$ is different of zero for some $P\in\Har{2j}(\RR^n)$.
    So, $m\le a$.

    Now, the theta series of the shadow is
    \begin{align*}
       \Theta_{\Sh(\Lambda),P} &= \Sh\Theta_{\Lambda,P}
                                = c_b(P)\,\Th2^{n+4j-8b}\Sh\De8^b+\cdots \\
                               &= 2^{n+4j-8b}(-1/16)^{b} c_b(P)\, q^{(n+4j-8b)/4} + \cdots,
    \end{align*}
    hence $2^{n+4j-8b} (-1/16)^{b} c_b(P)=\sum_{x\in\Sh(\Lambda)_{(n+4j-8b)/4}}P(x)$
    is different of zero for some $P\in\Har{2j}(\RR^n)$.
    So, $\sigma(\Lambda) = {n-8k} \le {n+4j-8b}$, and therefore $b\le k+j/2$.
  \end{proof}

  We give now similar statements for \emph{even} selfdual lattices.
  (The two propositions above are naturally true also for these lattices;
  though less precise.)
  We do not give the proofs, since the arguments are essentially the same as for
  the equivalent statements for general selfdual lattices.

  \begin{Proposition}\label{CalcThetaSeriesEven}
    Let $\Lambda\sseq\RR^n$ be an even selfdual lattice
    of rank $n=8N$.
    Then there exist $c_i\in\ZZ$ such that
    \begin{gather*}
      \Theta_\Lambda = \Q^N + c_1\De{24}\Q^{N-3} + \cdots +c_k \De{24}^{k}\Q^{N-3k},\quad k=[N/3],
    \end{gather*}
    with
    \[
      c_1 = n(h-30), \qquad \text{where } h:=\card{\Lambda_2}/n.
    \]
  \end{Proposition}

  \begin{Proposition}\label{CalcThetaSeriesEvenBis}
    Let $\Lambda\sseq\RR^n$ be an even selfdual lattice of of rank $n=8N$
    and of minimum~$m=2M$.
    \begin{enumerate}
      \item For every even positive integer~$j$,
          there exist linear forms $c_i:\Har{2j}(\RR^n)\to\CC$ such that
          \[
            \Theta_{\Lambda,P} = \sum_{i=M}^{\![(N+j/2)/3]\!} c_i(P)\,\De{24}^i\Q^{N+j/2-3i},
            \qquad\forall P\in\Har{2j}(\RR^n).
          \]
          In particular, if $3M>N+j/2$,
          then $\Theta_{\Lambda,P}=0$ for every $P\in\Har{2j}(\RR^n)$.
      \item For every odd positive integer~$j$,
          there exist linear forms $c_i:\Har{2j}(\RR^n)\to\CC$ such that
          \[
            \Theta_{\Lambda,P} = \sum_{i=M}^{\![(N+j/2)/3]\!} c_i(P)\,\R\,\De{24}^i\Q^{N+(j-3)/2-3i},
            \qquad\forall P\in\Har{2j}(\RR^n).
          \]
          In particular, if $3M>N+(j-3)/2$,
          then $\Theta_{\Lambda,P}=0$ for every $P\in\Har{2j}(\RR^n)$.
    \end{enumerate}
  \end{Proposition}

\section{Zero coefficients of modular forms}\label{SectZeroCoeff}

  From Lemma~\ref{CharactShellIsDesign},
  in order to find spherical design strengths of the shells of a lattice,
  we have to look for vanishing coefficients of $\Theta_{\Lambda,P}$
  for $P$~harmonic homogeneous polynomials of different degrees.
  We give here the results concerning
  vanishing coefficients of modular forms of the form $\Ph^{\epsilon}\Th3^{\alpha}\De8^{\beta}$,
  since we will meet them several times later. Recall that $\NN=\{0,1,2,\ldots\}$.

  \begin{Lemma}\label{ZeroCoeffModularForms}
    Among the coefficients of the modular forms
    \[
      \Ph^{\epsilon}\Th3^{\alpha}\De8^{\beta} = \sum_{m\in\alpha+\NN}a_m q^m
      \qquad \text{and}\qquad
      \Sh(\Ph^{\epsilon}\Th3^{\alpha}\De8^{\beta}) = \sum_{\!m\in\alpha/4+2\NN\!}a_m q^m,
    \]
    where $\epsilon\in\{0,1\}$, $\alpha\ge0$, $\beta\ge0$,
    the following are equal to zero:
    \begin{enumerate}
      \newcommand*\settabs{\kern70pt \= \kern 120pt \= \kern70pt \= \kill}
      \item[(a)]
        \begin{tabbing}\settabs
          $(\Th3^4\De8)^k$ : \> $m-k\equiv1\bmod2$,    \> $\Sh\bigl((\Th3^{16}\De8)^k\bigr)$ : \> $m\equiv2\bmod4$;
        \end{tabbing}
      \item[(b)]
        \begin{tabbing}\settabs
          $\Th3$              : \> $m\ne a^2$,         \> $\Sh(\Th3^2)$          : \> $m\ne(a^2+b^2)/4$,   \nl
          $\Sh(\Th3)$         : \> $m\ne a^2/4$,       \> $\Th3^3$               : \> $m=4^a(8b+7)$,       \nl
          $\Th3^2$            : \> $m\ne a^2+b^2$,     \> $\Ph\,\Th3^{12}$       : \> $m=1$;
        \end{tabbing}
      \item[(c)]
        \begin{tabbing}\settabs
          $\Th3  \De8$        : \> $m=4^a(8b+5)$,      \> $\Ph\,\Th3^3\De8$            : \> $m=4^a(8b+7)$, \nl
          $\Th3^2\De8$        : \> $m\ne a^2+b^2$,     \> $\Ph\,\Th3^{16}\De8$         : \> $m=4^a2$,      \nl
          $\Sh(\Th3^2\De8)$   : \> $m\ne(a^2+b^2)/4$,  \> $\Sh(\Ph\,\Th3^{16}\De8)$    : \> $4^a 2$,       \nl
          $\Th3^3\De8$        : \> $m=4^a(8b+7)$,      \> $\Sh(\Ph\,\Th3^{40}\De8)$    : \> $m=24$,        \nl          
          $\Th3^7\De8$        : \> $m=4^a(8b+3)$;    
        \end{tabbing}
      \item[(d)]
        \begin{tabbing}\settabs
          $\Th3^5\De8^2$           : \> $m=4^a(8b+1)$, \> $\Ph\,\Th3^{20}\De8^2$       : \> $m=3$,       \nl
          $\Th3^{12}\De8^2$        : \> $m=4^a$,       \> $\Ph\,\Th3^{33}\De8^2$       : \> $m=4$,       \nl
          $\Ph\,\Th3^8\De8^2$      : \> $m=4^a$,       \> $\Sh(\Ph\,\Th3^{33}\De8^2$)  : \> $m=49/4$,    \nl
          $\Sh(\Ph\,\Th3^8\De8^2$) : \> $m=4^a$;
        \end{tabbing}
      \item[(e)]
        \begin{tabbing}\settabs
          $\Th3^4\De8^3$           : \> $m=4^a2$,      \>  $\Sh(\Ph\,\Th3^{24}\De8^3)$ : \> $m=2^a$,       \nl
          $\Ph\,\Th3^{24}\De8^3$   : \> $m=2^a$.
        \end{tabbing}      
    \end{enumerate}
  \end{Lemma}

  \begin{Remark}
    We have checked numerically that,
    for modular forms of the Lemma with
    $\beta\le3$ and $\alpha\le36$,
    there is no other zero coefficient $a_m$
    for $m\le1200$.
  \end{Remark}

  \begin{proof}[Sketch of the proof]
    Let $\Theta(z)$ be one of the series of the Lemma.

    If $\Theta(z)$ is of the form
    $\Ph^{\epsilon}\Th3^{\alpha}\De8^{\beta}$
    or of the form $\Sh(\Ph^{\epsilon}\Th3^{4\gamma}\De8^{\beta})$,
    we can express $\Theta(z)$
    as a polynomial function of $\Th3(z)$ and $\Th4(z)$.
    If $\Theta(z)$ is of the form $\Sh(\Ph^{\epsilon}\Th3^{\alpha}\De8^{\beta})$,
    with $\alpha$ an integer which is not a multiple of~$4$,
    we can express $\Theta(4z)$ as a polynomial function of $\Th3(z)$ and $\Th4(z)$,
    by using the formulae of \cite[Chap.~4, \S4.1, p.~104]{ConwaySloane}.
    In the sequel of the proof, we suppose that we are in the first case;
    the second case is treated in a similar way, replacing $\Theta(z)$ by $\Theta(4z)$.
     
    Let be $\omega:=e^{i\pi/4}$.
    For $c$ an integer between $0$ and $7$, we define
    \begin{align*}
      \Theta_c(z) & := \frac18 \sum_{k=0}^7 \omega^{-ck}\,\Theta\biggl(\frac{z+k}{4}\biggr) \\
      &= \frac18 \sum_{k=0}^7 \omega^{-ck}\sum_{m\ge0}a_me^{i\pi m(z+k)/4} \\
      &= \frac18 \sum_{k=0}^7 \omega^{-ck}\sum_{j=0}^7 \omega^{jk}
          \sum_{\!\!\substack{m\ge0\\ m\equiv j\bmod8}\!\!} a_m e^{i\pi mz/4} \\
      &= \sum_{j=0}^7
         \biggl(
         \underbrace{\frac18\sum_{k=0}^7 \omega^{(j-c)k}}
               _{\substack{\text{$0$ if $j\ne c$}\\\text{$1$ if $j=c$}}}
         \biggr)
         \sum_{\!\!\substack{m\ge0\\ m\equiv j\bmod8}\!\!} a_me^{i\pi mz/4} \\
      &= \sum_{\!\!\substack{m\ge0\\ m\equiv c\bmod8}\!\!} a_me^{i\pi mz/4}.
    \end{align*}
    Therefore, whenever $\Theta_c(z)=0$, we have $a_m=0$
    for $m\equiv c\bmod 8$.
    Using
    the identities given in  \cite[Chap.~4, \S4.1, p.~104]{ConwaySloane},
    we can express
    \[
      \Theta_c(z) = \Th2(z)^d\,F\bigl(\Th3(z),\Th4(z)\bigr), \qquad 0\le d\le 3,\quad d\equiv c \bmod 8,
    \]
    where $F$ is a polynomial; so we can check whether $\Theta_c(z)=0$.
    (A computer software like Maple is highly recommended
    for performing these calculations.)

    We have 
    \begin{gather*}
      \Theta_{\text{even}}(z):=\frac12 \Bigl(\Theta(z) + \Theta(z+1) \Bigr)
          = \sum_{\!\!m\equiv 0\bmod2\!\!} \!a_{m}e^{i\pi mz},  \\[\jot]
      \Theta_{\text{odd}}(z):=\frac12 \Bigl(\Theta(z) - \Theta(z+1) \Bigr)
          = \sum_{\!\!m\equiv 1\bmod2\!\!} \!a_{m}e^{i\pi mz}, \\[\jot]
      \Theta_0(z) = \sum_{\!\!m\equiv 0\bmod2\!\!} \!a_{4m}e^{i\pi mz}, \qquad
      \Theta_4(z) = \sum_{\!\!m\equiv 1\bmod2\!\!} \!a_{4m}e^{i\pi mz},
    \end{gather*}
    and $\Theta_{\text{even}}$, $\Theta_{\text{odd}}$, $\Theta_0$ and $\Theta_4$
    can be expressed in function of $\Th3$ and $\Th4$.
    If $\Theta_{\text{even}}=0$ [respectively $\Theta_{\text{odd}}=0$],
    then $a_m=0$ for $m$ even [resp.~odd].
    If $\Theta_{\text{even}}$ [respectively $\Theta_{\text{odd}}$]
    is a multiple of $\Theta_0$ [resp.~$\Theta_4$], we have a relation
    between $a_m$ and $a_{4m}$ for $m$~even [resp.~odd],
    which allows to decide whenever $a_{4^k m}=0$.

    These considerations suffice for most of the theta series mentioned in the Lemma.
    For $\Th3^2\De8$ [resp.~$\Sh(\Th3^2\De8)$], we use Lemma~\ref{ThetaSeriesOfCubicLattice} below
    to show that $a_m=0$ whenever the shell of $\Z2$ [resp.~$\Sh(\Z{2})$] of norm~$m$ is empty.
  \end{proof}

  The following lemma provides another method to show that certain coefficients of modular forms are nonzero.
  
  \begin{Lemma}\label{LemmaPowerOfSeries}
    Let $\varphi^{(0)}$ and $\psi$ be two formal series in $q$ with integral coefficients such that
    \begin{gather*}
      \varphi^{(0)} =  a^{(0)}_k q^k + \sum_{j\ge k+1}a^{(0)}_j q^j, \qquad a^{(0)}_k\ge1,\ a^{(0)}_j\in\ZZ,  \\
      \psi  = b_0 + b_1 q + \sum_{j\ge2} b_j q^j,   \qquad b_0,b_1\ge1,\ b_j\ge0,
    \end{gather*}
    and let
    \[
      \varphi^{(n)}:=\varphi^{(0)}\,\psi^n = \sum_{j\ge k} a^{(n)}_j q^j.
    \]
    Then the sequence $M_n$, $n\ge0$ defined by
    \[
      M_n:=\max\{ m\in\NN \mid a^{(n)}_j>0,\ k\le j\le m \}
    \]
    is nondecreasing and unbounded.
  \end{Lemma}
  
  \begin{proof}
    We have $M_0\ge k$. For every $n\ge0$, we have $\varphi^{(n+1)}=\varphi^{(n)}\psi$; thus,
    for $k\le j\le M_n$,
    \[
      a^{(n+1)}_j = a^{(n)}_j b_0 +  \sum_{i\ge1} a^{(n)}_{j-i} b_i \ge a^{(n)}_j b_0 > 0.
    \]
    Consequently, $M_{n+1}\ge M_n$.
    Moreover, for $M=M_n$, we have
    \[
      a^{(n+1)}_{M+1} = a^{(n)}_{M+1} b_0 + a^{(n)}_M b_1 + \sum_{i\ge2} a^{(n)}_{M+1-i} b_i
               \ge a^{(n)}_{M+1} b_0 + a^{(n)}_M b_1 > a^{(n)}_{M+1}.
    \]
    Consequently, $a^{(n+1)}_{M+1}\ge a^{(n)}_{M+1}+1$, since all coefficients are integers.
    Similarly,
    \[
      a^{(n+h)}_{M+1}\ge a^{(n+h-1)}_{M+1}+1\ge\dots\ge a^{(n)}_{M+1}+h.
    \]
    Therefore, for $h$ large enough, we have $a^{(n+h)}_{M+1}>0$, and then $M_{n+h}\ge M_n+1$.
    That shows that the sequence $M_n$ is unbounded.
  \end{proof}    

\section{Root systems of norm 2}\label{SectRootSystems}

  For analysing selfdual lattices of rank at most~$24$,
  we need some informations on the first shells of a selfdual lattice~$\Lambda$.
  The decomposition $\Lambda=\Z{p}\oplus L$ gives us the exact form of the
  shell~$\Lambda_1$; we give now the form of~$\Lambda_2$.

  Let us first recall that a \emph{root system}
  is  a subset~$R\sseq\RR^n$ such that
  \begin{enumerate}
    \item $R$ is finite and does not contain~$0$;
    \item for every $x,y\in R$, we have $\refl{y}(x)\in R$,
          where $\refl{y}$ is the reflection of axis $\RR y$;
    \item for every $x,y\in R$, the number $2\scal{x}{y}/\scal{x}{x}$
          is an integer.
  \end{enumerate}

  {\bf Warning.} In the usual definition, 
  it is required that $R$ span $\RR^n$ as a vector space.
  Here, we leave out this condition in order to simplify the formulation of our results;
  in particular, we consider the empty set in $\RR^n$ as a root system.

  \begin{LemmaDefinition}\label{SecondLayerIsRootSystem}
    Let $\Lambda$ be an integral lattice. Then $\Lambda_2$ is a root system,
    called the \emph{root system of $\Lambda$}.
  \end{LemmaDefinition}

  \begin{proof}
    For $x$ and $y\in\Lambda_2$,
    we have $\refl{y}(x)= x -  \scal{x}{y} y$.
    Therefore, since $\scal xy$ is an integer, we have $\refl{y}(x)\in\Lambda_2$.
    Moreover, $2\scal{x}{y}/\scal{x}{x} = \scal{x}{y}$ is clearly an integer.
  \end{proof}

  If $\Lambda$ is a selfdual lattice of minimum at least~$2$ and of root system~$R$,
  we write sometimes $\Lambda=R^+$.
  (It happens that this notation is unambiguous up to rank~$23$,
  that is there is at most one selfdual lattice of minimum at least~$2$
  whose root system is a given root system in $\RR^{n}$, $n\le23$.)

  \medbreak
  Let us now recall some classical facts on root systems.
  If $R\sseq\RR^p$ and $S\sseq\RR^q$ are root systems, their \emph{orthogonal union}
  is the root system in $\RR^{p+q}$ defined by
  \[
    R+S := \bigl( R \oplus \{0\} \bigr) \cup \bigl( \{0\} \oplus S \bigr) \sseq \RR^p\oplus\RR^q.
  \]
  We write $k\,R$ for $R+\cdots+R$, $k$~terms.
  A root system
  is called \emph{irreducible} if it is not
  an orthogonal union of smaller root systems.
  Clearly, any root system can be written uniquely
  (up to permutation of the terms) as an orthogonal union of irreducible
  root systems.
  Note that a nonempty irreducible root system always spans its ambient space.
  The only empty irreducible root system is of dimension~1, and is noted $\OO_1$.

  An important number for an irreducible root system~$R$ is its \emph{Coxeter number},
  which is the integer~$h$ satisfying the relation
  \[
    \card{R} = nh.
  \]
  If $R$ is empty, we have $h=0$.

  The list of irreducible root systems is well known.
  We give here the list of those of \emph{norm 2},
  that is those root systems~$R$ that verify $\scal xx=2$ for every $x\in R$
  (they are also called \emph{simply laced} root systems by some authors).
  In the notation, the index indicates the dimension of the space where the root system lies.
  \begin{tabbing}
    \kern20pt \= $\D_n$, $n\ge4$,\quad \=  \kern120pt \= $\E_8$,\qquad \= \kill
    \> $\OO_1$,          \> $h=0$;      \> $\E_6$,          \> $h=12$;  \\
    \> $\A_n$, $n\ge1$,  \> $h=n+1$;    \> $\E_7$,          \> $h=18$;  \\
    \> $\D_n$, $n\ge4$,  \> $h=2(n-1)$; \> $\E_8$,          \> $h=30$.
  \end{tabbing}
  We denote by $\OO_n:=n\,\OO_1$ the empty root system in $\RR^n$.

  The following definition, justified by the next lemma,
  is inspired by
  the corresponding notion for lattices
  (see~\cite[pp.~28ff]{MartinetVenkov}).

  \begin{Definition}\label{DefStronglyEutacticRootSystem}
        A root system~$R$ is called \emph{strongly eutactic} if,
        in its decomposition in irreducible root systems,
        all irreducible components have the same Coxeter number.
        In this case, we define the \emph{Coxeter number} of~$R$ as the Coxeter number
        of any of its irreducible components.
  \end{Definition}
  A nonempty strongly eutactic root system spans its ambiant space.  
  Note that the equality $\card{R} = nh$ holds for strongly eutactic root systems.

  \begin{Lemma}
    Let $R$ be a root system of norm~$2$
    that is a spherical $3$-design.
    Then $R$ is strongly eutactic.
  \end{Lemma}

  The converse is true also: see~Proposition~\ref{ClassificationRootSystems}.i.

  \begin{proof}
    Let us write a root system $R\sseq\RR^n$ as
    \[
      R = R_1 + R_2 + \cdots + R_k,
    \]
    where the $R_i$'s are irreducible.
    Let $V_i$ be the subspace of $\RR^n$ where $R_i$ lies
    (thus $\RR^n=V_1\oplus V_2\oplus\cdots\oplus V_k$),
    and let $n_i=\dim V_i$.
    Any point $x\in\RR^n$, is written uniquely as
    \[
      x = x_1 + x_2 + \cdots + x_k, \qquad x_i\in V_i.
    \]
    Let us consider the harmonic polynomials of degree~$2$ defined by
    \[
      f_{i,j}(x) = \frac{1}{2n_i}\scal{x_i}{x_i} - \frac{1}{2n_j}\scal{x_j}{x_j}.
    \]
    We have
    \[
      \sum_{x\in R} f_{i,j}(x) = \frac{2\card{R_i}}{2n_i} - \frac{2\card{R_j}}{2n_j} = h_i-h_j.
    \]
    So, if $R$ is a spherical $3$-design, we must have $\sum_{x\in R} f_{i,j}(x)=0$
    and therefore $h_i=h_j$; in other words, $R$ is strongly eutactic.
  \end{proof}

  We recall now the notion of reproducing kernel,
  which will help us to analyse strongly eutactic root systems:
  
  Let $\mathcal{H}$ be a complex (or a real) finite-dimensional Hilbert space of functions
  on a nonempty set~$\Omega$.
  We use the convention that hermitian scalar products are antilinear in the \emph{first} variable.
  There exists a unique function $\Phi:\Omega\times\Omega\to\CC$,
  called \emph{reproducing kernel}, such that
  $\Phi(x,\cdot)\in\mathcal{H}$ for all $x\in\Omega$, and
  \[
    f(x) = \scal{\Phi(x,\cdot)}{f}, \qquad \forall f\in\mathcal{H},\ \forall x\in\Omega.
  \]
  This kernel verifies $\overline{\Phi(y,x)}=\Phi(x,y)$.
  It is of \emph{positive type}; that is, for any finitely supported function
  $\Omega\to\CC$, $y\mapsto\lambda_y$,
  \[
    \sum_{x,y\in\Omega} \overline{\lambda_x}\lambda_y \Phi(x,y) \ge0.
  \]
  Moreover, the set $\{\Phi(x,\cdot)\mid x\in\Omega\}$ generates $\mathcal{H}$.
  For all this, see for example \cite{BekkaHarpe}.

  The positivity of $\Phi$ implies the following result:

  \begin{Lemma}\label{CharactPositivityOfRK}
    Let $\Omega\to\CC$, $y\mapsto \lambda_y$ be a finitely supported function. Then
    \[
      \sum_{x,y\in\Omega} \overline{\lambda_x}\,\lambda_y\,\Phi(x,y) =0
      \qquad\text{if and only if}\qquad
      \sum_{y\in\Omega} \lambda_y\,f(y)=0,\quad\forall f\in\mathcal{H}.
    \]
  \end{Lemma}

  \begin{proof}
     Let $\Omega'$ be the vector space of finitely supported functions on $\Omega$.
     We define a hermitian form~$h$ on $\Omega'$ by:
    \[
      h\bigl(\mu,\lambda\bigr)
              = \sum_{x,y\in\Omega} \overline{\mu_x}\,\lambda_y\,\Phi(x,y),
      \qquad \mu,\lambda\in\Omega'.
    \]
    For $x\in \Omega$, let $\delta_x\in\Omega'$ denote the function on $\Omega$
    that takes value~$1$ at~$x$ and value~$0$ elsewhere.
    Let $\Omega_0'$ be the set of $\lambda\in\Omega'$ such that
    $\scal{\mu}{\lambda}=0$ for every $\mu\in\Omega'$.
    Then $\lambda\in\Omega'_0$ if and only if $h(\delta_x,\lambda)=0$
    for all $x\in X$, if and only if
    \begin{gather*}
      \sum_{y\in\Omega} \lambda_y \Phi(x,y) = 0, \quad\forall x\in\Omega.
    \end{gather*}
    But since $\{\Phi(x,\cdot)\mid x\in \Omega\}$ generates $\mathcal{H}$,
    the last condition is equivalent to
    \[
      \sum_{y\in\Omega} \lambda_y f(y) = 0, \quad\forall f\in\mathcal{H}.
    \]
    Now, the positivity of $\Phi$ implies that
    \[
      h(\lambda,\lambda)=0 
      \qquad\text{if and only if}\qquad
      \lambda \in \Omega_0'.
    \]
    This is exactly what is claimed by the Lemma.
  \end{proof}

  \medbreak

  Let us now consider the following special case:
  \begin{gather*}
    \Omega = \Se,\qquad \mathcal{H}=\Har{j}(\RR^n),\\
    \scal{f}{g}=\int_{\Se} \overline{f(u)}\,g(u)\,d\sigma(u),\quad f,g\in\Har{j}(\RR^n),
  \end{gather*}
  where $\sigma$ is the probability measure on $\Se$ invariant by rotation.
  Let $\Phi^{(j)}$ be the corresponding reproducing kernel.
  It is known that 
  \begin{equation}\tag{$*$}
    \Phi^{(j)}(x,y) = \Qp{j}\bigl(\scal{x}{y}\bigr) \qquad \forall x,y\in\Se,
  \end{equation}
  where $\Qp{j}(t)$ is an appropriate Gegenbauer polynomials,
  with the normalisation
  $\Qp{j}(1)=\dim\bigl(\Har{j}(\RR^n)\bigr)$;
  see \cite{DelsarteGoethalsSeidel} and \cite[\S IX.3]{Vilenkin}.
  We have, for example,
  \begin{gather*}
    \Qp{0}(t) = 1, \qquad
    \Qp{1}(t) = nt, \qquad
    \Qp{2}(t) = \frac{n+2}{2}\bigl(nt^2 - 1\bigr), \\
    \Qp{3}(t) = \frac{n(n+4)}{6}  \Bigl( (n+2)t^3 - 3t \Bigr), \\
    \Qp{4}(t) = \frac{n(n+6)}{24} \Bigl( (n+2)(n+4)t^4 - 6(n+2)t + 3 \Bigr).
  \end{gather*}

  We can now prove:

  \begin{Proposition}\label{ClassificationRootSystems}
    Let $R$ be a nonempty strongly eutactic root system of norm~$2$.
    Let us consider the conditions
    \begin{equation*}\tag{$C_{2j}$}
      \sum_{x\in R} f(x) = 0, \qquad \forall f\in\Har{2j}(\RR^n).
    \end{equation*}
    Then:
    \begin{enumerate}
      \item
        Condition $(C_2)$ always holds
        (equivalently, nonempty strongly eutactic root systems are spherical $3$-designs);
      \item
        Condition $(C_4)$ holds if and only if
        $R$ is equivalent to one of the following systems:
        \[
           \A_1,\quad \A_2,\quad \D_4,\quad \E_6,\quad \E_7,\quad \E_8;
        \]
      \item
        Condition $(C_6)$ holds if and only if
        $R$ is equivalent to one of the following systems:
        \[
           \A_1,\quad 2\A_1,\quad \E_8,\quad 2\E_8,\quad \D_{16};
        \]
      \item
        Condition $(C_8)$ holds if and only if
        $R$ is equivalent to one of the following systems:
        \[
          \A_1,\quad\A_2;
        \]
      \item
        Condition $(C_{10})$ holds if and only if
        $R$ is equivalent to one of the following systems:
        \[
          \A_1,\quad2\A_1,\quad\A_2\quad\D_4,\quad \E_8;
        \]
      \item
        Condition $(C_{12})$ holds if and only if
        $R$ is equivalent to the system
        \[
          \A_1.
        \]
    \end{enumerate}
  \end{Proposition}

  \begin{proof}
    For $x$ and $y\in R$, we have $\scal{x}{y}\in\{0,\pm1,\pm2\}$.
    Let $x\in R$, and let, for $\alpha\in\{0,\pm1,\pm2\}$,
    \[
      N_\alpha := \card{ \{y\in\R \mid \scal{x}{y}=\alpha\} }.
    \]
    We have the evident relations $N_{-\alpha}=N_\alpha$, $N_2=1$,
    and $N_0+2N_1+2N_2=nh$. Moreover,
    it is known that $N_\alpha$ is independent of~$x$ and that
    \begin{equation*}
      N_1 = 2h-4
    \end{equation*}
    (see  Bourbaki~\cite{BourbakiLie}, Lie~VI, \S~1.11, prop.~3.2,
    where it is stated for (nonempty) irreducible root systems;
    but it immediately extends to nonempty strongly eutactic root systems).

    Now, we renormalize $R$ by $\widetilde{R}=\frac{1}{\sqrt2}R$, so that $\widetilde{R}\sseq\Se$.
    By Lemma~\ref{CharactPositivityOfRK}
    (applied with $\lambda_x=1$ if $x\in X$ and $\lambda_x=0$ otherwise)
    and equation~$(*)$,
    Condition~$(C_{2j})$ is equivalent to
    \[
      \sum_{x,y\in\widetilde{R}}\Qp{2j}\bigl(\scal{x}{y}\bigr) = 0.
    \]
    However, since for every $x\in R$ we have [note that $\Qp{2j}(-\alpha)=\Qp{2j}(\alpha)$]:
    \begin{align*}
       \sum_{y\in\widetilde{R}}\Qp{2j}\bigl(\scal{x}{y}\bigr)
         &= \sum_{\alpha}N_{\alpha}\Qp{2j}(\alpha/2) \\
         &= 2\Qp{2j}(1) + (4h-8)\Qp{2j}(1/2) + (nh-4h+6)\Qp{2j}(0),
    \end{align*}
    Condition~$(C_{2j})$ is equivalent to
    \[
      2\Qp{2j}(1) + (4h-8)\Qp{2j}(1/2) + (nh-4h+6)\Qp{2j}(0) =0.
    \]
    This condition is linear in~$h$ and polynomial in~$n$
    (since $\Qp{2j}$ is polynomial in~$n$).
    It is, for $2j\le12$:
    \begingroup
      \renewcommand*\nl{\displaybreak[0]\\[2\jot]}
      \begin{gather*}
        0 = 0  \tag{$C_2$} \nl
        n(n+4)(n+6)\bigl((n-10)h+6(n+2)\bigr) = 0 \tag{$C_4$} \nl
        n(n+2)(n+6)(n+10)\bigl((n^2-48n+272)h+30(n-4)(n+4)\bigr) = 0 \tag{$C_6$} \nl
        \begin{split}
          & n(n+2)(n+4)(n+8)(n+14) \\
          & \quad\times \bigl( (n-4)(n-30)(n-50)h + 42(n+6)(3n^2-14n+40)\bigr) = 0
        \end{split}
        \tag{$C_8$} \nl
        \begin{split}
          & (n-2)n(n+2)(n+4)(n+6)(n+10)(n+18) \\
          & \quad\times \bigl( (n-24)(n-28)(n-76)h + 30(n+8)(17n^2-8n+336)\bigr) = 0
        \end{split}
        \tag{$C_{10}$} \nl
        \begin{split}
          & n(n+2)(n+4)(n+6)(n+8)(n+12)(n+22) \\
          & \quad\times \bigl( (n^5-186n^4+10852n^3-228504n^2+1659232n-967680)h \\[-1pt]
          & \qquad\qquad\quad+ 66n(n-2)(n+10)(31n^2+130n+1144)\bigr) = 0
        \end{split}
        \tag{$C_{12}$}
      \end{gather*}
    \endgroup
    We can now explicit the positive integral solutions of these equations.
    Let us take for example~$(C_4)$: we have $({n-10})h+6({n+2})=0$, hence
    \[
      h=\frac{6(n+2)}{10-n}.
    \]
    Since $h>0$, we have $n<10$. Now, if we introduce $n=1,2,\ldots,9$ in the last formula,
    we obtain the integral solutions
    \[
      (n,h) = (1,2)\,,\ (2,3)\,,\ (4,6)\,,\ (6,12)\,,\ (7,18)\,,\ (8,30)\,,\ (9,16).
    \]
    It is now easy to list all strongly eutactic root systems having one of these parameters.
  \end{proof}

\part{Application to some selfdual lattices}

\section{The cubic lattices}\label{SectCubicLattices}

  Let
  \[
    \Z{n} := \{x=(x_1,\ldots,x_n)\in\RR^n \mid x_i\in\ZZ,\ i=1,\cdots,n\}
  \]
  be the \emph{cubic lattice} of rank~$n$. It is an odd selfdual lattice whose shadow is:
  \[
    \Sh(\Z{n}) = \{x\in\RR^n \mid x_i\in\ZZ+{\textstyle\frac12},\ i=1,\cdots,n\}.
  \]
  In particular, $\sigma(\Z{n})=n$.


  \begin{Lemma}\label{ThetaSeriesOfCubicLattice} We have
    \[
      \Theta_{\Z{n},P} =
        \begin{cases}
           \Th3^n                  & \text{if $P=1$,}              \nl
           0                       & \text{if $P\in\Har2(\RR^n)$,} \nl
           c_1(P)\,\De8\Th3^n      & \text{if $P\in\Har4(\RR^n)$,} \nl
           c_2(P)\,\Ph\,\De8\Th3^n & \text{if $P\in\Har6(\RR^n)$,}
        \end{cases}
    \]
    where $c_1$ is a nonzero linear form if and only if $n\ge2$,
    and $c_2$ is a nonzero linear form if and only if $n\ge3$.
  \end{Lemma}%

  \begin{Remark}
    The theta series of $\Sh(\Z{n})$ are given by Proposition~\ref{ThetaSeriesOfShadow}.
  \end{Remark}

  \begin{proof}
    The equality $\Theta_{\Z{n}}=\Th3^n$ follows directly from
    $\Theta_{\Z{n}}=(\Theta_{\Z{}})^n$ and $\Theta_{\Z{}}=\Th3$.
    Proposition~\ref{CalcThetaSeriesBis} gives the form of $\Theta_{\Z{n},P}$ for $P\in\Har{2j}(\RR^n)$,
    $j=1,2,3$.
    
    It remains to show that
    $c_1$ and $c_2$ are not identically zero for $n$ large enough,
    i.e., there exists
    $P\in\Har4(\RR^n)$ (if $n\ge2$) such that $\Theta_{\Z{n},P}\ne0$,
    and there exists $Q\in\Har6(\RR^n$) (if $n\ge3$) such that $\Theta_{\Z{n},Q}\ne0$.
    We can chose for example
    \[
      P(x) = x_1^4 + x_2^4 - 6 x_1^2 x_2^2, \qquad \Theta_{\Z{n},P}(z)= 4q+\cdots,
    \]
    and
    \begin{gather*}
       Q(x)=     (x_1^6+x_2^6+x_3^6)
              - 15(x_1^4x_2^2+x_2^4x_3^2+x_3^4x_1^2)
              + 90(x_1^2x_2^2x_3^2), \\
      \Theta_{\Z{n},Q}(z)= 6q+\cdots.
    \end{gather*}

    Finally, the fact that $c_1\equiv0$ [respectively $c_2\equiv0$] if $n=1$ [respectively if $n=1$~or~$2$]
    is an easy exercise.
  \end{proof}
  

  The zero coefficients of the Fourier expansion of the series of the last Lemma
  are given in Lemma~\ref{ZeroCoeffModularForms}.
  Therefore, by Lemma \ref{CharactShellIsDesign},
  we have:

  \begin{Theorem}\label{StrengthCubicLattices}\
    \NOPAGEBREAK
    \let\settabs=\setstandardtabs
    \begin{enumerate}
      \item
        For $n\ge2$, all nonempty shells of $\Z{n}$ and $\Sh(\Z{n})$
        are spherical $3$-designs.
      \item
        The following shells are spherical $3\half$-designs:\NOPAGEBREAK
        \begin{tabbing}\settabs
           \> $\Z{2}$ and $\Sh(\Z{2})$ : \> every nonempty shell, \\
           \> $(\Z{16})_m$             : \> $m=4^a 2$, $a\ge0$    \\
           \> $\Sh(\Z{16})_m$          : \> $m=4^a 2$, $a\ge1$    \\
           \> $\Sh(\Z{40})_m$          : \> $m=24$. 
        \end{tabbing}
      \item
        The following shells are spherical $5$-designs:\NOPAGEBREAK
        \begin{tabbing}\settabs
           \> $(\Z{4})_m$     : \> $m=2a$,        $a\ge1$   \\
           \> $(\Z{7})_m$     : \> $m=4^a(8b+3)$, $a,b\ge0$ \\
           \> $\Sh(\Z{16})_m$ : \> $m=4a+2$,      $a\ge1$.
        \end{tabbing}
    \end{enumerate}
  \end{Theorem}

  Note that the zero coefficients of series $\De8\Th3^n$ for $n=2,3$,
  and $\Ph\,\De8\Th3^n$ for $n=3$,
  that are mentioned in Lemma~\ref{ZeroCoeffModularForms},
  correspond to empty shells of the lattices $\Z{2}$ and $\Z{3}$.
  Zero coefficients of $\De8\Th3^n$ for $n=1$ and
  $\Ph\,\De8\Th3^n$ for $n=1,2$ are irrelevant for our purpose,
  because they correspond to cases where $c_1$ and $c_2$
  are zero in Lemma~\ref{ThetaSeriesOfCubicLattice}.
  Similar remarks hold for the shadows of these modular forms.
  
  We give in Appendix an alternative proof that some shells of the cubic lattices of rank 4~and~7
  are $5$-designs.
  
  \begin{Proposition}\label{NonStrengthCubicLattices}\ \NOPAGEBREAK
    \begin{enumerate}
    \item
      For $n\ge2$ and $1\le m\le 1200$,
      the nonempty shells of norm~$m$ of $\Z{n}$
      are not spherical $5$-designs,
      except those mentioned in Theorem~\ref{StrengthCubicLattices}.
    \item
      For $n\ge2$ and $n/4\le m\le n/4+1200$,
      the nonempty shells of norm~$m$ of $\Sh(\Z{n})$ are not a spherical $5$-designs,
      except those mentioned in Theorem~\ref{StrengthCubicLattices}.
   \end{enumerate}
  \end{Proposition}
  
  \begin{proof}
   (i) We check first numerically that the statement holds for $n<408$.
       Then, we apply Lemma~\ref{LemmaPowerOfSeries} with $\varphi^{(0)}=\De8(z)=q+O(q^2)$
       and $\psi=\Th3(z)=1+2q+O(q^4)$. Numerical computations give $M_{408}\ge1200$, and 
       thus the coefficients $a^{(n)}_j$ of
       \[
         \varphi^{(n)} = \De8(z)\,\Th3(z)^n = \sum_{j\ge1}a^{(n)}_j q^j
       \]
       are positive for $1\le j\le 1200$ and $n\ge 408$.
       By Lemma~\ref{ThetaSeriesOfCubicLattice} and Lemma~\ref{CharactShellIsDesign},
       the corresponding shells $(\Z{n})_j$ are not spherical $5$-designs.
       
  (ii) We check first numerically that the statement holds for $n<426$.
      Then, we apply Lemma~\ref{LemmaPowerOfSeries} to
      \begin{gather*}
        \varphi^{(0)}= -16\,\Sh\De8(z/2) = 1 -16q + O(q^2), \\
        \psi=q^{-1/8}\Th2(z/2) = 2 + 2q + O(q^3).
      \end{gather*}
      (The substitution $z\mapsto z/2$ is equivalent to the substitution $q\mapsto q^{1/2}$.)
      Numerical computations give $M_{426}\ge600$.
      It follows that the coefficients $a^{(n)}_j$ of
      \[
         \varphi^{(n)}=-16\,q^{n/8}\,\Sh\De8(z/2)\,\Th2(z/2)^n = \sum_{j\ge0}a^{(n)}_j q^j \\
      \]
      are positive for $0\le j\le 600$ and $n\ge 426$. In other words, the coefficients $c^{(n)}_j$ of
      \[
         \Sh\De8(z)\,\Th2(z)^n = \sum_{j\ge0}c^{(n)}_j q^{2j+n/4} \\
      \]
      are negative for $n/4\le 2j+n/4 \le n/4+1200$ and $n\ge426$.
      The proof is achieved by invoking Lemma~\ref{ThetaSeriesOfCubicLattice}
      together with Proposition~\ref{ThetaSeriesOfShadow},
      and Lemma~\ref{CharactShellIsDesign}.
  \end{proof}

\section{The Witt lattices}\label{SectWittLattices}

  Let $n$ be a positive multiple of~$4$. The \emph{Witt lattice} of rank~$n$
  is the lattice
  \[
    \W{n} := \Biggl\{x=(x_1,\dots,x_n)\in\RR^n \Biggm|
                                    \begin{aligned}
                                       &2x_i\in\ZZ\quad \forall i,   \\[-\jot]
                                       &x_i-x_1\in\ZZ\quad \forall i,\\[-\jot]
                                       &x_1+x_2+\dots+x_n\in 2\ZZ
                                     \end{aligned}
                                     \Biggr\}.
  \]
  It is selfdual, and
  it is even if $n\equiv0\bmod8$ and odd otherwise.
  When $n\equiv4\bmod8$, the shadow of $\W{n}$ is
  \[
    \Sh(\W{n}) = \Biggl\{x=(x_1,\dots,x_n)\in\RR^n \Biggm|
                                    \begin{aligned}
                                       &2x_i\in\ZZ\quad \forall i,   \\[-\jot]
                                       &x_i-x_1\in\ZZ\quad \forall i,\\[-\jot]
                                       &x_1+x_2+\dots+x_n\in 2\ZZ+1
                                     \end{aligned}
                                     \Biggr\}.
  \]
  The lattice $\W{4}$ is equivalent to the cubic lattice $\Z{4}$ and is
  thus analysed in the previous Section; therefore,
  we assume that $n\ge8$.
  The lattice $\W{8}$ is analysed more precisely in the next Section.

  \begin{Lemma} Let $n$ be a multiple of~$4$ greater than or equal to~$8$. Then
    \[
      \Theta_{\Lambda,P} = \begin{cases}
                              \frac12(\Th2^n+\Th3^n+\Th4^n) & \text{if $P=1$,} \\[\jot]
                              0                             & \text{if $P\in\Har2(\RR^n)$,} \\[\jot]
                              c_1(P)\,\Th2^4\Th3^4\Th4^4(-\Th2 ^{n-4}+\Th3^{n-4}-\Th4^{n-4})
                                                            & \text{if $P\in\Har4(\RR^n)$,} \\[\jot]
                              c_2(P)\,\Th2^4\Th3^4\Th4^4\Bigl(   (\Th3^4+\Th4^4)\Th2 ^{n-4}  \\
                                                     \qquad   + (\Th4^4-\Th2^4)\Th3 ^{n-4}
                                                              - (\Th2^4+\Th3^4)\Th4 ^{n-4}
                                                      \Bigr)
                                                            & \text{if $P\in\Har6(\RR^n)$.}
                           \end{cases}
    \]
    where $c_1$ and $c_2$ are nonzero linear forms on $\Har4(\RR^n)$ and $\Har6(\RR^n)$ respectively. 
  \end{Lemma}

  \begin{proof}
    Let $G$ the subgroup of $\Orth(n)$ containing the transformations of the form
    \[
      (x_1,\dots,x_n) \mapsto (\epsilon_1x_1,\dots,\epsilon_n x_n)
    \]
    where $\epsilon_i=\pm1$ and $\epsilon_1\epsilon_2\dotsm\epsilon_n=1$.
    It is clear that $G$ leaves $\W{n}$ invariant.
    We denote by $\Har{2j}(\RR^n)^G$ the set of elements of
    $\Har{2j}(\RR^n)$ which are invariant under the action of~$G$
    given by $(\gamma P)(x)=P(\gamma^{-1}x)$.
    Let
    \[
      \pi : \Har{2j}(\RR^n) \to  \Har{2j}(\RR^n)^G,
      \qquad
      \pi(P) = \frac{1}{\card{G}} \sum_{\gamma\in G} \gamma P
    \]
    be the projection on the invariant part of $\Har{2j}(\RR^n)$.
    Since for all $\gamma\in G$
    \[
      \Theta_{\W{n},\gamma P} = \Theta_{\gamma^{-1}\W{n},P} = \Theta_{\W{n},P},
    \]
    we have $\Theta_{\W{n},P} = \Theta_{\W{n},\pi(P)}$.
    Therefore, it suffices to prove the Lemma for $P\in\Har{2j}(\RR^n)^G$.

    Let
    \[
      \varW{n} = \Biggl\{x=(x_1,\dots,x_n)\in\RR^n \Biggm|
                                      \begin{aligned}
                                         &2x_i\in\ZZ\quad \forall i,   \\[-\jot]
                                         &x_i-x_1\in\ZZ\quad \forall i,\\[-\jot]
                                         &-x_1+x_2+\dots+x_n\in 2\ZZ
                                       \end{aligned}
                                       \Biggr\},
    \]
    which is a lattice equivalent to $\W{n}$.
    Let $\D_n:=\{x\in\Z{n}\mid\scal xx\in 2\ZZ\}$;
    we have, by definition, $\Sh(\Z{n})=\D_n^\sharp\setminus\D_n$.
    We have the following inclusions of lattices,
    where labels indicate indices of sublattices:
    $$
      \xymatrix{   & \D_n^\sharp \ar@{-}[dl]_2 \ar@{-}[d]^2  \ar@{-}[dr]^2  \\
                \W{n} & \Z{n} & \varW{n} \\
                   & \D_n \ar@{-}[ul]^2 \ar@{-}[u]_2  \ar@{-}[ur]_2  \\
               }
    $$
    Expressing $\D_n^\sharp$ as the union of its classes modulo $\D_n$, we obtain
    \[
      \D_n^\sharp = (\W{n}\setminus\D_n)\sqcup(\Z{n}\setminus\D_n)\sqcup(\varW{n}\setminus\D_n)\sqcup\D_n.
    \]
    Hence, since $\Sh(\Z{n}) = \D_n^\sharp\setminus\Z{n}$,
    \[
      \Sh(\Z{n})  = (\W{n}\setminus\D_n)\sqcup(\varW{n}\setminus\D_n),
    \]
    from which we deduce,
    \[
      \Theta_{\Sh(\Z{n})  ,P} = \Theta_{\W{n},P} + \Theta_{\varW{n},P} - 2\Theta_{\D_n,P}.
    \]
    As $2\Theta_{\D_n,P}(z)=\Theta_{\Z{n},P}(z)+\Theta_{\Z{n},P}(z+1)$, we have
    \[
     \Theta_{\W{n},P}(z) + \Theta_{\varW{n},P}(z)
       = \Theta_{\Z{n},P}(z)+\Theta_{\Z{n},P}(z+1) + \Theta_{\Sh(\Z{n}),P}(z).
    \]
    Let $\sigma$ be the orthogonal transformation of $\RR^n$ defined by
    \[
      \sigma(x_1,x_2,\dots,x_n)=(-x_1,x_2,\dots,x_n);
    \]
    we have $\sigma\varW{n}=\W{n}$.
    Let $P\in\Har{2j}(\RR^n)^G$, where $0\le 2j\le 6$.
    It is easily checked that, since $2j\le n$,
    $P$ is a polynomial which is \emph{even} in $x_1$, $x_2$, \dots,~$x_n$,
    and therefore $P$ is invariant under the action of $\sigma$.
    Thus, $\Theta_{\varW{n},P}=\Theta_{\W{n},\sigma P}=\Theta_{\W{n},P}$, and    
    \[
     \Theta_{\W{n},P}(z)
       = \frac12\Bigl(\Theta_{\Z{n},P}(z)+\Theta_{\Z{n},P}(z+1) + \Theta_{\Sh(\Z{n}),P}(z)\Bigr).
    \]
    We use Lemma~\ref{ThetaSeriesOfCubicLattice} and identities
    found in \cite[Chap.~4, \S4.1, p.~104]{ConwaySloane} to conclude.
  \end{proof}

  Now, using Lemma \ref{CharactShellIsDesign}, we have:

  \begin{Theorem}\label{StrengthWittLattices}
    Let $n$ be a multiple of~$4$ greater than or equal to~$8$.
    \NOPAGEBREAK
    \let\settabs=\setstandardtabs
    \begin{enumerate}
      \item
        All nonempty shells of $\W{n}$ and $\Sh(\W{n})$
        are spherical $3$-designs.
      \item
        All nonempty shells of $\W{8}$
        are spherical $7$-designs.
      \item
        All nonempty shells of $\W{16}$
        are spherical $3\half$-designs.
    \end{enumerate}
  \end{Theorem}
  
  We do not know if any shell of $\W{n}$, for $n\ge12$ is a spherical $5$-design;
  however, as in the case of cubic lattices, we can show:

  \begin{Proposition}\label{NonStrengthWittLattices}
    For $n\ge12$ and $1\le m\le1200$,
    the nonempty shells of norm~$m$ of
    $\W{n}$ and of their shadows are not
    spherical $5$-designs.
  \end{Proposition}

  \begin{proof}
    For $m< n/4$, we remark that
   $(\W{n})_m=(\Z{n})_m$ when $m$ is an even integer and $(\W{n})_m=\emptyset$ otherwise.
   If moreover $n\equiv4\bmod 8$, then we have $\bigl(\Sh(\W{n})\bigr)_m=(\Z{n})_m$ when $m$ is an odd integer
   and $\bigl(\Sh(\W{n})\bigr)_m=\emptyset$ otherwise.
   Therefore, for $n>4800$, it suffices to apply Proposition~\ref{NonStrengthCubicLattices}.
   For $n\le 4800$, we content ourself with a numerical verification.
  \end{proof}

  The case of $\W{8}$ is considered in more details in the next Section.

\section{Even selfdual lattices of rank at most 24}\label{SectEvenLattices24}

  We recall here the classification of these lattices, due to Niemeier.
  Note that we already know by Proposition~\ref{BasicFactsShadow}
  that their rank is a multiple of~$8$.
  
  \begin{Theorem}\label{ClassificationEvenSelfdualLattices}\
    \NOPAGEBREAK
    \begin{enumerate}
      \item There is exactly one even selfdual lattice of rank~$8$, that is
        the Korkine-Zolotarev lattice
        \[
          \W{8} = \E_8^+.
        \]
      \item There are exactly two even selfdual lattices of rank~$16$, that is
        \[
          \W{8}\oplus\W{8} = (2\,\E_8)^+ \qquad\text{and}\qquad \W{16} = \D_{16}^+ .
        \]
      \item There is a bijection between the even selfdual lattices of rank~$24$
        and
        the twenty-four strongly eutactic root systems of norm~$2$ and of rank~$24$,
        given by $\Lambda\leftrightarrow\Lambda_2$.
        These root systems are (in parentheses is the Coxeter number of the system):
        \begin{gather*}
          \OO_{24}              \ ( 0), \qquad
          24\,\A_1              \ ( 2), \qquad
          12\,\A_2              \ ( 3), \qquad
          8\,\A_3               \ ( 4), \\
          6\,\A_4               \ ( 5), \qquad
          4\,\A_5 + \D_4        \ ( 6), \qquad
          6\,\D_4               \ ( 6), \qquad
          4\,\A_6               \ ( 7), \\
          2\,\A_7 + 2\,\D_5     \ ( 8), \qquad
          3\,\A_8               \ ( 9), \qquad
          9\,\A_9 + \D_6        \ (10), \qquad
          4\,\D_6               \ (10), \\
          \A_{11} + \D_7 + \E_6 \ (12), \qquad
          4\,\E_6               \ (12), \qquad
          2\,\A_{12}            \ (13), \qquad
          3\,\D_8               \ (14), \\
          \A_{15} + \D_9        \ (16), \qquad
          \D_{10} + 2\,\E_7     \ (18), \qquad
          \A_{17} + \E_7        \ (18), \qquad
          2\,\D_{12}            \ (22), \\
          \A_{24}               \ (25), \qquad
          \D_{16}+\E_8          \ (30), \qquad
          3\,\E_8               \ (30), \qquad
          \D_{24}               \ (46).  
        \end{gather*}
    \end{enumerate}
  \end{Theorem}

  Claims (i)~and~(ii) are consequences of Claim~(iii), since
  $\Lambda=R^+$ is an even unimodular lattice of rank~$8$ [respectively $16$]
  if and only if $\Lambda\oplus\W{8}\oplus\W{8} = (R+2\,\E_8)^+$
  [resp.\ $\Lambda\oplus\W{8} = (R+\E_8)^+$] is an even unimodular
  lattice of rank~$24$.
  
  The lattice corresponding to $\OO_{24}$ is the famous Leech lattice,
  and those corresponding to nonempty strongly eutactic root systems of rank~$24$
  are the the Niemeier lattices;
  those corresponding to the three last root systems listed are
  $\W{16}\oplus\W{8}$, $\W{8}\oplus\W{8}\oplus\W{8}$, and $\W{24}$ respectively.

  There exist several proofs of Claim~(iii);
  a proof explaining the bijection is found in \cite{Venkov24}.

  \begin{Lemma}\label{ThetaSeriesOfEvenSelfdualLattices}\
    \NOPAGEBREAK
    \begin{enumerate}
      \item
        Let $\Lambda=\W{8}=\E_8^+$ be the Korkine-Zolotareff lattice.
        Then
        \[
          \Theta_{\Lambda,P} =
            \begin{cases}
              \Q                      & \text{if $P=1$,} \\
              0                       & \text{if $P\in\Har{2j}(\RR^8)$, $2j=2,4,6,10$,} \\
              c_1(P)\,\De{24}         & \text{if $P\in\Har{8}(\RR^8)$,} \\
              c_2(P)\,\Q\,\De{24}     & \text{if $P\in\Har{12}(\RR^8)$,} \\
              c_3(P)\,\R\,\De{24}     & \text{if $P\in\Har{14}(\RR^8)$,} \\
              c_4(P)\,\Q^2\De{24}     & \text{if $P\in\Har{16}(\RR^8)$,} \\
              c_5(P)\,\R\,\Q^2\De{24} & \text{if $P\in\Har{18}(\RR^8)$,}
            \end{cases}
        \]
        where $c_1$, $c_2$, $c_3$, $c_4$, and  $c_5$
        are nonzero linear form on
        $\Har{8}(\RR^8)$, $\Har{12}(\RR^8)$, $\Har{14}(\RR^8)$, $\Har{16}(\RR^8)$, and $\Har{18}(\RR^8)$
        respectively.
      \item
        Let $\Lambda$ be one of the two even selfdual lattices of rank~$16$.
        Then
        \[
          \Theta_{\Lambda,P} =
            \begin{cases}
              \Q^2                    & \text{if $P=1$,} \\
              0                       & \text{if $P\in\Har{2j}(\RR^{16})$, $2j=2,6$,} \\
              c_1(P)\,\De{24}         & \text{if $P\in\Har{4}(\RR^{16})$,} \\
              c_2(P)\,\Q\,\De{24}     & \text{if $P\in\Har{8}(\RR^{16})$,} \\
              c_3(P)\,\R\,\De{24}     & \text{if $P\in\Har{10}(\RR^{16})$,} \\
              c_4(P)\,\Q^2\De{24}     & \text{if $P\in\Har{12}(\RR^{16})$,} \\
              c_5(P)\,\R\,\Q^2\De{24} & \text{if $P\in\Har{14}(\RR^{16})$,}
            \end{cases}
        \]
        where $c_1$, $c_2$, $c_3$, $c_4$, and  $c_5$
        are nonzero linear forms on
        $\Har{4}(\RR^{16})$, $\Har{8}(\RR^{16})$, $\Har{10}(\RR^{16})$, $\Har{12}(\RR^{16})$, and $\Har{14}(\RR^{16})$
        respectively.
      \item
        Let $\Lambda$ be the Leech lattice.
        Then
        \[
          \Theta_{\Lambda,P} =
            \begin{cases}
              \Q^3- 720\De{24}         & \text{if $P=1$,} \\
              0                        & \text{if $P\in\Har{2j}(\RR^{24})$, $2j=2,4,6,8,10,14$,} \\
              c_1(P)\,\De{24}^2        & \text{if $P\in\Har{12}(\RR^{24})$,} \\
              c_2(P)\,\Q\,\De{24}^2    & \text{if $P\in\Har{16}(\RR^{24})$,} \\
              c_3(P)\,\R\,\De{24}^2    & \text{if $P\in\Har{18}(\RR^{24})$,} \\
              c_4(P)\,\Q^2\De{24}^2    & \text{if $P\in\Har{20}(\RR^{24})$,} \\
              c_5(P)\,\R\,\Q^2\De{24}  & \text{if $P\in\Har{22}(\RR^{24})$,}
            \end{cases}
        \]
        where $c_1$, $c_2$, $c_3$, $c_4$, and  $c_5$
        are nonzero linear forms on the corresponding spaces.
      \item
        Let $\Lambda$ be an even selfdual lattice of rank~$24$
        and of minimum~$2$.
        Let $h:=\card{\Lambda_2}/24$.
        Then
        \[
          \Theta_{\Lambda,P} =
            \begin{cases}
              \Q^3- 24(30-h)\De{24}    & \text{if $P=1$,} \\
              0                        & \text{if $P\in\Har{2}(\RR^n)$,} \\
              c_1(P)\,\Q\,\De{24}      & \text{if $P\in\Har{4}(\RR^n)$,} \\
              c_2(P)\,\R\,\De{24}      & \text{if $P\in\Har{6}(\RR^n)$,} \\
              c_3(P)\,\Q^2\De{24}      & \text{if $P\in\Har{8}(\RR^n)$,} \\
              c_4(P)\,\R\,\Q^2\De{24}  & \text{if $P\in\Har{10}(\RR^n)$,}
            \end{cases}
        \]
        where $c_1$, $c_2$, $c_3$, and $c_4$ are nonzero linear forms on
        the corresponding spaces.
    \end{enumerate}
  \end{Lemma}

  \begin{proof}
    It suffices to apply Propositions \ref{CalcThetaSeriesEven}~and~\ref{CalcThetaSeriesEvenBis}.
    In order to see that the linear forms $c_i$, $i=1,2,\ldots$, are nonzero,
    we apply Lemma~\ref{CharactShellIsDesign}
    and Proposition~\ref{ClassificationRootSystems} to $\Lambda_2$ if $\Lambda$ is not the Leech lattice,
    and to $\Lambda_4$ if $\Lambda$ is the Leech lattice.
    (Note that, since $\{\scal{x}{y}\mid x,y\in\Lambda_4\} = \{0,\pm1,\pm2,\pm4\}$,
    the shell $\Lambda_4$ is a \emph{tight} $11$-design: see
    \cite[Theorems 5.12~and~6.8]{DelsarteGoethalsSeidel}
    or \cite[Theorems 5.3~and~5.4]{GoethalsSeidel}.
    It follows that $\Lambda_4$ cannot be a $12$-design.)
  \end{proof}

  \begin{Remark}
    In fact, Propositions \ref{CalcThetaSeriesEven}~and~\ref{CalcThetaSeriesEvenBis}
    show that the root system of an even
    selfdual lattice of rank~$8$ [respectively $16$, $24$]
    is nonempty and a $7\half$-design [resp. is nonempty and a $5\half$-design, is strongly eutactic].
    Theorem~\ref{ClassificationEvenSelfdualLattices} says that,
    for each root system in dimension up to~$24$ which satisfies the above condition,
    there exists exactly one even selfdual lattice with such a root system.
    But in higher dimensions, there can be more than one even selfdual lattice with the same root system.
  \end{Remark}
  
  As in the preceding case, Lemma~\ref{CharactShellIsDesign} give:\NOPAGEBREAK

  \begin{Theorem}\label{StrengthEvenUnimodularLattices}\
    \NOPAGEBREAK
    \begin{enumerate}
      \item
        All shells of $\W{8}=\E_8^+$
        are spherical $7\half$-designs.
      \item
        All shells of the two even selfdual lattices of rank~$16$
        are spherical $3\half$-designs.
      \item
        All shells of the Leech lattice are spherical $11\half$-designs.
      \item
        All shells of any even selfdual lattice of rank~$24$ and of minimum~$2$
        are spherical $3$-designs.
    \end{enumerate}
  \end{Theorem}
  
  \begin{Remark}
    We have checked numerically that the shells of norm at most~$1200$
    of these lattices
    are not spherical designs of higher strength. 
  \end{Remark}
  
  Claims (i) to (iii) of the Theorem above are special cases
  of Theorem~16.4 of \cite{VenkovMartinet1} (see also \cite{VenkovExtremal}),
  since the lattices mentioned in these claims are \emph{extremal}.
  See the remark at the end of Section~\ref{SectClassifModularForms}.

  Another consequence of our analysis is a reformulation
  in terms of spherical design strength of shells of lattices
  of a famous
  conjecture of Lehmer, which states that
  the Ramanujan numbers $\tau(m)$ are never zero for $m\ge1$.
  This conjecture has been verified for $m\le 10^{15}$
  \cite[\S~3.3]{Serre85}.

  \begin{Proposition}\label{TheoremOnLehmerConjecture}
        Let $\tau(m)$, $m\ge1$, be the Ramanujan numbers defined by
        \[
          \De{24}(z)=\sum_{m\ge1} \tau(m)q^{2m}.
        \]
        For $m\ge1$, the following are equivalent:\NOPAGEBREAK
        \begin{enumerate}
          \item[(a)] $\tau(m)=0$;
          \item[(b)] $(\W{8})_{2m}$ is an $8$-design (and therefore an $11$-design);
          \item[(c)] for any even selfdual lattice~$\Lambda$ of rank~$16$,
                     $\Lambda_{2m}$ is a $4$-design (and therefore a $7$-design);
        \end{enumerate}
  \end{Proposition}

  Therefore, Lehmer's conjecture is true if and only if no shell of $\W{8}$
  is an $8$-design,
  and if and only if no shell of any even selfdual lattice of rank~$16$
  is a $4$-design.
  
  Similar conjectures could be stated for other modular forms than $\De{24}$,
  and have an equivalent formulation in terms of spherical design strength of shells of lattices
  (see Proposition~C in the Introduction).

\section{Selfdual lattices with long shadow}\label{SectLatticesLongShadow}

In this section, we consider selfdual lattices with $\sigma(\Lambda)=n-\nobreak8$.
We begin with those of minimum at least~$2$.
We recall here their classification \cite{Elkies8}.

\begin{Theorem}\label{ClassificationSelfdualLatticesType8}
  There is a bijection between the selfdual lattices with $\sigma(\Lambda)=n-\nobreak8$
  of minimum at least~$2$
  and the 14 strongly eutactic root systems satisfying $h=2(23-n)$,
  given by $\Lambda\leftrightarrow\Lambda_2$.
  These root systems are (in parentheses are the corresponding values of $n$~and~$h$):
  \begin{gather*}
    \E_8           \ ( 8,30),\qquad
    \D_{12}        \ (12,26),\qquad
    2\,\E_7        \ (14,18),\qquad
    \A_{15}        \ (15,16),\\
    2\,\D_8        \ (16,14),\qquad
    \A_{11} + \E_6 \ (17,12),\qquad
    2\,\A_9        \ (18,10),\\
    3\,\D_6        \ (18,10),\qquad
    2\,\A_7 + \D_5 \ (19, 8),\qquad
    4\,\A_5        \ (20, 6),\\
    5\,\D_4        \ (20, 6),\qquad
    7\,\A_3        \ (21, 4),\qquad
    22\,\A_1       \ (22, 2),\qquad
    \OO_{23}       \ (23, 0).
  \end{gather*}
\end{Theorem}

Recall that we denote by $R^+$ the (unique up to dimension~23) selfdual lattice~$\Lambda$
of root system $R=\Lambda_2$.
Apart from $\W{8}=\E_8^+$, which is even, all these lattices are odd. They are of minimum~$2$,
except $\OO_{23}^+$ which is of minimum~$3$ (it is the so-called \emph{shorter Leech lattice}).

The theta series of such lattices are now easy to calculate
(we exclude here the theta series of $\W{8}=\E_8^+$, which have been given in
Lemma~\ref{ThetaSeriesOfEvenSelfdualLattices}).

\begin{Lemma}\label{ThetaSeriesOfLatticesType8}
  Let $\Lambda$ be an odd selfdual lattice with $\sigma(\Lambda)=n-\nobreak8$
  and of minimum at least~$2$.
  \begin{enumerate}
    \item If $\Lambda$ is of minimum~$2$, then
      \[
        \Theta_{\Lambda,P} = 
          \begin{cases}
                  \Th3^n - 2n\,\Th3^{n-8}\De8   & \text{if $P=1$,} \\
                  0                             & \text{if $P\in\Har{2}(\RR^n)$,} \\
                  c_1(P)\,\Th3^{n-8}\De8^2      & \text{if $P\in\Har{4}(\RR^n)$,} \\
                  c_2(P)\,\Ph\Th3^{n-8}\De8^2   & \text{if $P\in\Har{6}(\RR^n)$.} \\
           \end{cases}
      \]
      where $c_1$ and $c_2$ are nonzero linear forms on $\Har{4}(\RR^n)$ and $\Har{6}(\RR^n)$ respectively. 
    \item If $\Lambda=\OO_{23}^+$ is the shorter Leech lattice, then 
      \[
        \Theta_{\Lambda,P} = 
          \begin{cases}
                  \Th3^{23} - 46\,\Th3^{15}\De8  & \text{if $P=1$,} \\
                  0                              & \text{if $P\in\Har{2j}(\RR^{23})$, $2j=2,4,6$,} \\
                  c_1(P)\,\Th3^{15}\De8^3        & \text{if $P\in\Har{8}(\RR^{23})$,} \\
                  c_2(P)\,\Ph\Th3^{15}\De8^3     & \text{if $P\in\Har{10}(\RR^{23})$,} \\
           \end{cases}
      \]
      where $c_1$ and $c_2$ are nonzero linear forms on $\Har{8}(\RR^n)$ and $\Har{10}(\RR^n)$ respectively. 
  \end{enumerate}
\end{Lemma}

\begin{proof}
  Propositions \ref{CalcThetaSeries} and~\ref{CalcThetaSeriesBis} give the forms of the theta series
  (incidently, they show that $\Lambda_2$ is strongly eutactic and that $h=2(23-n)$
  as stated in Theorem~\ref{ClassificationSelfdualLatticesType8}),
  except the fact that the linear forms $c_1$ and $c_2$ are nonzero.
  
  If $\Lambda$ is of minimum~$2$, by Proposition~\ref{CharactShellIsDesign} applied to the root system of $\Lambda$,
  the linear forms $c_1$ and $c_2$ are nonzero.
  
  Consider now the case $\Lambda=\OO_{23}^+$.
  Let $x\in\Lambda_3$, and set $N_\alpha:=\Card{\{y\in\Lambda_3 \mid \scal{x}{y}=\alpha\}}$.
  It is known that
  \[
    N_3=N_{-3}=2, \qquad
    N_1=N_{-1}=891, \qquad
    N_0=2816,
  \]
  and $N_\alpha=0$ for other~$\alpha$'s. (This follows from he fact that $\Lambda_3$ is a \emph{tight} $7$-design:
  see
  \cite[Theorems 5.12~and~6.8]{DelsarteGoethalsSeidel}
  or \cite[Theorems 5.3~and~5.4]{GoethalsSeidel}.)
  From this, using Lemma~\ref{CharactPositivityOfRK}
  (with $\lambda_x=1$ if $x\in\Lambda_3$, and $\lambda_x=0$ otherwise),
  it is easy to see that $c_1$ and $c_2$ are indeed nonzero.
\end{proof}

The zero coefficients of the Fourier series for the modular forms of the Lemma
are given by Lemma~\ref{ZeroCoeffModularForms}. We deduce:

\begin{Theorem}\label{StrengthLongShadowLattices}\
  \NOPAGEBREAK
  \let\settabs=\setstandardtabs
  \begin{enumerate}
    \item Every nonempty shell of any selfdual lattice with $\sigma(\Lambda)=n-\nobreak8$
       of minimum~$2$,
       and every nonempty shell of its shadow,
       is a spherical $3$-design.
    \item The following shells are spherical $3\half$-designs:
      \begin{tabbing}\settabs
        \> $(2\,\D_8)^+$                 : \> $m=4^a$,\ $a\ge1$
      \end{tabbing}
    \item The following shells are spherical $5$-designs:
      \begin{tabbing}\settabs
        \> $(4\,\A_5)^+$, $(5\,\D_4)^+$   : \> $m=4^a$,\  $a\ge1$ \\[\jot]
        \> $(2\,\D_8)^+$                 : \> $m=2a+1$,\ $a\ge1$
      \end{tabbing}
    \item All nonempty shells of the shorter Leech lattice $\OO_{23}^+$
      and of its shadow are spherical $7$-designs.
  \end{enumerate}
\end{Theorem}

\begin{Remark}
  We have checked numerically that the shells of norm at most~$1200$
  of these lattices and of their shadows
  are not spherical designs of higher strength. 
\end{Remark}

Let us now turn to selfdual lattices of minimum~$1$.

\begin{Lemma}
  Let $\Lambda$
  be a selfdual lattice with $\sigma(\Lambda)=n-\nobreak8$ of minimum~$1$.
  \begin{enumerate}
    \item If $\Lambda=\Z{1}\oplus\E_8^+ = \Z{1}\oplus\W{8}$, then
      \begin{equation*}
        \Theta_{\Lambda,P} =
          \begin{cases}
            \Th3^9 - 16\,\Th3\De8                    & \text{if $P=1$,} \\
            c_1(P)\,\Ph\Th3\De8                        & \text{if $P\in\Har2(\RR^9)$,} \\
            c_2(P)\,(\Th3^9\De8 + 8\,\Th3\De8^2)      & \text{if $P\in\Har4(\RR^9)$,} \\
            c_4(P)\,(\Ph\Th3^{9}\De8 -\Ph\Th3\De8^2) & \text{if $P\in\Har6(\RR^9)$,}
          \end{cases}
      \end{equation*}
      where $c_1$, $c_2$ and $c_4$ are nonzero linear forms on $\Har2(\RR^9)$,
      $\Har4(\RR^{9})$ and $\Har6(\RR^{9})$ respectively.
    \item If $\Lambda=\Z{1}\oplus\OO_{23}^+$, then
      \begin{equation*}
        \Theta_{\Lambda,P} =
          \begin{cases}
            \Th3^{24} - 46\,\Th3^{16}\De8                      & \text{if $P=1$,} \\
            c_1(P)\,\Ph\Th3^{16}\De8                             & \text{if $P\in\Har2(\RR^{24})$,} \\
            c_2(P)\,(\Th3^{24}\De8 -40\,\Th3^{16}\De8^2)        & \text{if $P\in\Har4(\RR^{24})$,} \\
            c_4(P)\,(\Ph\Th3^{24}\De8 -16\,\Ph\Th3^{16}\De8^2) & \text{if $P\in\Har6(\RR^{24})$,}
          \end{cases}
      \end{equation*}
      where $c_1$, $c_2$ and $c_4$ are nonzero linear forms on $\Har2(\RR^{24})$,
      $\Har4(\RR^{24})$ and $\Har6(\RR^{24})$ respectively.
    \item Otherwise, we have
      \begin{equation*}
        \Theta_{\Lambda,P} =
          \begin{cases}
            \Th3^n - 2n\,\Th3^{n-8}\De8                      & \text{if $P=1$,} \\
            c_1(P)\,\Ph\Th3^{n-8}\De8                           & \text{if $P\in\Har2(\RR^n)$,} \\
            c_2(P)\,\Th3^{n}\De8 + c_3(P)\,\Th3^{n-8}\De8^2    & \text{if $P\in\Har4(\RR^n)$,} \\
            c_4(P)\,\Ph\Th3^{n}\De8 + c_5(P)\,\Ph\Th3^{n-8}\De8^2
                                                              & \text{if $P\in\Har6(\RR^n)$,}
          \end{cases}
      \end{equation*}
      where $c_1$ is a nonzero linear form on $\Har2(\RR^n)$,
      $c_2$ and $c_3$ are linearly independent linear forms on $\Har4(\RR^n)$,
      and $c_4$ and $c_5$ are linearly independent linear forms on $\Har6(\RR^n)$.
  \end{enumerate}
\end{Lemma}

\begin{proof}
  Let $\Lambda=\Z{p}\oplus L$, where $L$ is a lattice of minimum~$2$ and of rank~$N$,
  and let $h$ be the Coxeter number of the strongly eutactic root lattice $L_2$.
  Let $V_1$ [respectively $V_2$] be the space generated by $\Z{p}$ [resp.~$L$],
  so that $\RR^n=V_1\oplus V_2$.
  For $x\in\RR^n$, let $x_i\in V_i$ ($i=1,2$)
  such that $x=x_1+x_2$. Let $\omega_i(x):=\scal{x_i}{x_i}$.

  First, Proposition~\ref{CalcThetaSeries} gives the exact form of $\Theta_{\Lambda}$.

  Now, by Proposition~\ref{CalcThetaSeriesBis}, we have
  $\Theta_{\Lambda,P}=c_1(P)\,\Ph\Th3^{n-8}\De8$ for $P\in\Har2(\RR^n)$.
  The polynomial
  \[
    P:=\frac{1}{2p}\omega_1-\frac{1}{2N}\omega_2\in\Har2(\RR^n)
  \]
  verifies $\Theta_{\Lambda,P}(z)=q+O(q^2)\ne0$; thus $c_1$ is not identically equal to zero.

  Then, we have
  $\Theta_{\Lambda,P}=c_2(P)\,\Th3^{n}\De8 + c_3(P)\,\Th3^{n-8}\De8^2$
  for $P\in\Har4(\RR^n)$.
  Let be the polynomial
  \[
    Q:=\frac{1}{2p}\omega_1^2 - \frac{p+2}{pN}\omega_1\omega_2 + \frac{p+2}{2N(N+2)}\omega_2^2\in\Har4(\RR^n).
  \]
  We have
  \begin{align*}
    \Theta_{\Lambda,Q}(z)
      &= q + 4\Bigl((p-1) + \frac{h(p+2)}{2(N+2)}\Bigr)q^2 + O(q^3) \\
      &= \biggl(q + (2n-8) q^2 + O(q^3) \biggr)
          +  \biggl(\Bigl(2p-2N+4 + \frac{2h(p+2)}{N+2}\Bigr)q^2 + O(q^3)\biggr) \\
      &= \Th3^n(z)\,\De8(z)
          + \Bigl(2p-2N+4 + \frac{2h(p+2)}{N+2}\Bigr) \Th3^{n-8}(z)\,\De8^2(z).
  \end{align*}
  If $L$ is neither $\OO_{23}^+$ nor $\E_8^+$,
  then, by Lemma~\ref{ThetaSeriesOfLatticesType8},
  there exists a $R\in\Har4(V_2)\sseq\Har4(\RR^n)$
  such that $\Theta_{L,R}\ne 0$. So we have
  \[
    \Theta_{\Lambda,R}(z) = \Theta_{\Z{p},1}(z)\,\Theta_{L,R}(z)
    = c_4(R)\,q^2 + O(q^3) = c_4(R)\,\Th3^{n-8}(z)\,\De8^2(z),
  \]
  with $c_4(R)\ne 0$; that shows that $c_3$ and $c_4$ are linearly independant.
  
  If $\Lambda=\Z{p}\oplus\E_8^+$ and if $p\ge2$, we can
  find a $S\in\Har4(V_1)\sseq\Har4(\RR^n)$ such that $\Theta_{\Z{p},S}=\Th3^p\De8$. So we have
  \[
    \Theta_{\Z{p}\oplus\E_8^+,S} = \Theta_{\Z{p},S\vphantom{E_8^+}}\Theta_{\E_8^+,1}
    = \Th3^{p}\De8 (\Th3^8-16\,\De8) = \Th3^n\De8 - 16\,\Th3^{n-8}\De8.
  \]
  Therefore, if we compare
  \begin{gather*}
    \Theta_{\Z{p}\oplus\E_8^+,Q}= \Th3^n(z)\,\De8(z) + 8p\, \Th3^{n-8}(z)\,\De8^2(z).
    , \\
    \Theta_{\Z{p}\oplus\E_8^+,S} = \Th3^n\De8 - 16\,\Th3^{n-8}\De8,
  \end{gather*}
  we find that $c_2$ and $c_3$ are linearly independent.

  If $\Lambda=\Z{1}\oplus\E_8^+$, it can be shown that $c_2$ and $c_3$ are not linearly independent.
  We do not prove it, since it has no effect on the conclusions of the following Theorem.

  If $\Lambda=\Z{p}\oplus\OO_{23}^+$, we proceed as for $\Z{p}\oplus\E_8^+$.

  A similar method is used for computing $\Theta_{\Lambda,P}$ for $P\in\Har6(\RR^n)$.
  Note however that in the cases $\Z{2}\oplus\E_8^+$ and $\Z{2}\oplus\E_8^+$,
  the polynomial corresponding to $S$ is a polynomial of the form
  \[
    S' = \bigl(N\,\omega_1-(p+8)\,\omega_2\bigr)\,f \in\Har6(\RR^n),\qquad \text{where $f\in\Har4(V_1)$}.
  \]
  and such that $\Theta_{\Z{2},f}\ne0$, since $\Theta_{\Z{2},P}=0$ for every $P\in\Har6(V_1)$.
\end{proof}

Again, Lemma~\ref{ZeroCoeffModularForms} describes indices of vanishing coefficients for the series of the previous Lemma.
Therefore, we have:

\begin{Theorem}\label{StrengthLongShadowLatticesBis}\
  \NOPAGEBREAK
  \begin{enumerate}
    \item
      Let $\Lambda$ be a selfdual lattice of rank~$24$ with $\sigma(\Lambda)=24-\nobreak8=16$.
      Then the shells $\Lambda_m$ and $\Sh(\Lambda)_m$
      are spherical $3$-designs for $m=4^a2$, $a\ge0$
      (except the shells $\Sh(\Lambda)_2$ and  $(\Z{1}\oplus\OO_{23}^+)_2$, which are empty).
    \item
      Let $\Lambda=\Z{3}\oplus\W{8}$.
      Then the shells $\Lambda_m$ are spherical $3$-designs for $m=4^a(8b+7)$, $a,b\ge0$.
  \end{enumerate}
\end{Theorem}

\begin{Remark}
  We have checked numerically that the shells of norm at most~$1200$
  of these lattices and of their shadows
  are not spherical designs of higher strength. 
\end{Remark}

\section{Odd selfdual lattices of rank~24 and of minimum~at least~2}\label{SectOddLattices24}

We need to know the precise form of the shell of norm~$2$
of these lattices and their shadows:

\begin{Proposition}\label{OddSelfdualLatticesDimension24}
  Let $\Lambda$ be an odd selfdual lattice of rank~$24$ and of minimum~$2$.
  Then there exist a root system~$R$,
  and strongly eutactic root systems $S$ and $T$ of Coxeter number $h_S$ and $h_T$ respectively,
  such that
  \begin{myitemize}
    \item $R=S\cap T$,
    \item $ h_S + h_T = 3 h_R+2 $  where $h_R:=\card{R}/24$,
    \item $s\in S\setminus R,\ t\in T\setminus R \implies \scal{s}{t}=\pm\frac{1}{2}$,
    \item $\Lambda_2=R$, and $\Sh(\Lambda)_2=(S\setminus R) \cup (T\setminus R)$.
  \end{myitemize}
  Moreover, $T=R$ or $S=R$ if and only if $R$ is strongly eutactic.
\end{Proposition}

\begin{proof}[Sketch of the proof]
  Let $\Lambda$ be an odd selfdual lattice of rank~$24$ and of minimum~$2$.
  Let $\Lambda_0$ be the even sublattice of index~$2$,
  and let $\Lambda_0^\sharp$ be its dual.
  Then $\Lambda_0^\sharp/\Lambda_0\simeq\ZZ/2\ZZ\oplus\ZZ/2\ZZ$,
  and we have the following diagramm (integers indicate indices  of sublattices):
    $$
      
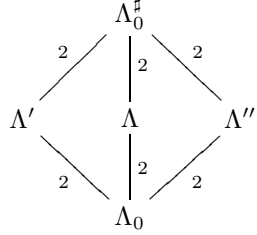
\xymatrix{   & \Lambda_0^\sharp \ar@{-}[dl]_2 \ar@{-}[d]^2  \ar@{-}[dr]^2  \\
                \Lambda' & \Lambda & \Lambda'' \\
                   & \Lambda_0 \ar@{-}[ul]^2 \ar@{-}[u]_2  \ar@{-}[ur]_2  \\
               }
    $$
  where $\Lambda'$ and $\Lambda''$ are even selfdual lattices. We have
  $
    \Sh(\Lambda) = (\Lambda'\setminus\Lambda_0) \cup (\Lambda''\setminus\Lambda_0).
  $
  The root systems of the Proposition are
  \[
    R = \Lambda_2, \qquad S=\Lambda'_2, \qquad T=\Lambda''_2,
  \]
  and, according to Theorem~\ref{ClassificationEvenSelfdualLattices}.iii,
  $S$ and $T$ are strongly eutactic.
  By Proposition~\ref{CalcThetaSeries}, the theta series of $\Lambda$ is
  \[
    \Theta_\Lambda = \Th3^n + c_1\De8\Th3^{n-8} + c_2 \De8^2\Th3^{n-16},
  \]
  where
  \begin{align*}
    c_2  & = 24(\card{\Lambda_2}/24 - 2) = 24(h_R-2) \\
         & = \card{\Sh(\Lambda)_{2}} = 24(h_S+h_T-2h_R),
  \end{align*}
  hence $h_S+h_T=3h_R+2$.

  For $x\in\RR^{24}$ and $\alpha\in\RR$,
  let $N_\alpha^{x,A}:=\card{\{y\in A\mid\scal{x}{y}=\alpha\}}$.
  Let $s\in S\setminus R$; it can be shown that (see Lemma~\ref{MutalScalarProductsOfRST} below):
  \begin{gather*}\tag{$*$}
    N_1^{s,R} = N_{-1}^{s,R} = N_1^{s,T} =  N_{-1}^{s,T}  = 3h_R - h_T = h_S-2, \\
    N_{1/2}^{s,T} = N_{-1/2}^{s,T} = 12(h_T-h_R), \\
    N_{0}^{s,R} = N_{0}^{s,T} = 8h_T + 6h_S - 12.
  \end{gather*}
  Let us suppose that $R$ is strongly eutactic and that $S\setminus R\ne\emptyset$.
  Let $s\in S\setminus R$;
  since the polynomial function  $x\mapsto\scal{s}{x}^2 - 2\scal{x}{x}$ is harmonic,
  we have the equality
  \[
    \sum_{r\in R}\scal{s}{r}^2 = 4h_R,
  \]
  hence
  \[
    2h_S-4 = N_1^{s,R} +N_{-1}^{s,R} = \sum_{r\in R}\scal{s}{r}^2 = 4h_R.
  \]
  Using the equality $h_S+h_T=3h_R+2$, we deduce that $h_T=h_R$, therefore $T=R$
  and $h_S=2h_R+2$.
  Conversely, if $T=R$, then $R$ is strongly eutactic.  
\end{proof}

To complete the proof, we establish formulae~$(*)$ above:

\begin{Lemma}\label{MutalScalarProductsOfRST}
  In the proof of the preceding Lemma, formulae $(*)$ hold.
\end{Lemma}

\begin{proof}[Sketch of the proof]
  Let $s\in S\setminus R$.
  The scalar product on $\Lambda_0^\sharp$ reduces to a bilinear form
  on $\Lambda_0^\sharp/\Lambda_0$ with values in $\bigl(\frac{1}{2}\ZZ\bigr)\big/\ZZ$,
  which is easy to explicit. In particular, for $t\in T$,
  we have:
  \[
    \scal{s}{t} = \begin{cases} \pm1/2 & \text{if $t\in T\setminus R$,} \\
                                0,\pm1 & \text{if $t\in R$.}
                  \end{cases}
  \]
  (The values $\pm3/2$ are excluded in the first case by
  computing the norm of $s\mp t$.)
  From this, we deduce $N^{s,T}_{\pm1/2}=12(h_T-h_R)$, and $N^{s,R}_\epsilon=N^{s,T}_\epsilon$
  for $\epsilon=0,\pm1$. Finally, $N^{s,T}_\epsilon$ is deduced from the strongly eutaxy of $T$.
  For example, we have
  \[    
    -\frac14 N^{s,T}_{-1/2} + \frac34 N^{s,T}_{1/2} + 2 N^{s,T}_1
     = \sum_{t\in T} \bigr(\scal{s}{t}+1\bigl)\scal{s}{t} = 2h_T.
  \]
  The last equality follows from the fact that $T$ is a spherical $2$-design
  (use Theorem~3.2 and Formula~3.6 in \cite{VenkovMartinet1}, or Proposition~1.13 in \cite{Gaeta}).
  Therefore $N^{s,T}_1=3h_R-h_T$.
\end{proof}

\medbreak

There are 156~odd selfdual lattices of rank~$24$, classified in~\cite{BorcherdsOdd24Ori}
and listed in \cite{BorcherdsOdd24}.
We give here the list of those with strongly eutactic root system,
because it is the occasion to point out a remarkable bijection between two sets
(compare with Theorems \ref{ClassificationEvenSelfdualLattices}.iii~and~\ref{ClassificationSelfdualLatticesType8}).

\begin{Theorem}
  There is a bijection between the set of
  odd selfdual lattices~$\Lambda$ of rank~$24$, of minimum~$2$,
  and of strongly eutactic root system,
  and the set of pairs of embedded strongly eutactic root systems $R\sseq S$ in $\RR^{24}$
  satisfying
  $h_S=2 h_R+2$.
  The bijection is given by $R=\Lambda_2$ and $S=\Sh(\Lambda)_2$.
  The list of pairs $R\sseq S$ is given below.
  (The numbers in brackets refer to Table~17.1 of \cite{BorcherdsOdd24}\/%
  \footnote{There are some errors in that table in the two first editions,
            which have been corrected  in the third edition.}
  and Table~III of \cite{Bacher}.
  These lattices correspond to bold edges in the neighbourhood graph of Figure~17.1 in \cite{BorcherdsOdd24}.)
  \begin{gather*}
    \OO_{24}       \ \sseq\  24\,\A_1           \ [1],    \qquad
    24\,\A_1       \ \sseq\  6\,\D_4            \ [6],    \qquad
    8\,\A_3        \ \sseq\  4\,\D_6            \ [32],   \\
    6\,\D_4        \ \sseq\  3\,\D_8            \ [74],   \qquad
    2\,\A_7 + \D_5 \ \sseq\  \D_{10} + 2\,\E_7  \ [105],  \qquad
    4\,\D_6        \ \sseq\  2\,\D_{12}         \ [130],  \\
    3\,\D_8        \ \sseq\  3\,\E_8            \ [145],  \qquad
    3\,\D_8        \ \sseq\  \E_8 + \D_{16}     \ [146],  \qquad
    2\,\D_{12}     \ \sseq\  \D_{24}            \ [154].
  \end{gather*}
\end{Theorem}

\medbreak

We are now ready to give the theta series for selfdual lattices of rank~$24$:\NOPAGEBREAK

\begin{Lemma}
  Let $\Lambda$ be an odd selfdual lattice of rank~$24$, and of minimum at least~$2$.\NOPAGEBREAK
  \begin{enumerate}
    \item If $\Lambda_2$ is not strongly eutactic, and if $h:=\card{\Lambda_2}/24$, then
      \begin{equation*}
        \Theta_{\Lambda,P} =
          \begin{cases}
            \Th3^{24} - 48\,\Th3^{16}\De8 + 24(h+2)\,\Th3^{8}\De8^2
                                                & \text{if $P=1$,} \\
            c_1(P)\,\Ph\Th3^8\De8^2             & \text{if $P\in\Har2(\RR^{24})$,} \\
            c_2(P)\,\Th3^{16}\,\De8^2 + c_3(P)\,\Th3^{8}\De8^3
                                                & \text{if $P\in\Har4(\RR^{24})$,}
          \end{cases}
      \end{equation*}
      where $c_1$ is a nonzero linear form on $\Har2(\RR^{24})$, and
      $c_2$ and $c_3$ are linearly independent linear forms on $\Har4(\RR^{24})$.
    \item If $\Lambda_2$ is nonempty and strongly eutactic of Coxeter number $h$, then
      \begin{equation*}
        \Theta_{\Lambda,P} =
          \begin{cases}
            \Th3^{24} - 48\,\Th3^{16}\De8 + 24(h+2)\,\Th3^{8}\De8^2
                                                & \text{if $P=1$,} \\
            0                                   & \text{if $P\in\Har2(\RR^{24})$,} \\
            c_2(P)\,\Th3^{16}\,\De8^2 + c_3(P)\,\Th3^{8}\De8^3
                                                & \text{if $P\in\Har4(\RR^{24})$,} \\
            c_4(P)\,\Ph\Th3^{16}\,\De8^2 + c_5(P)\,\Ph\Th3^{8}\De8^3
                                                & \text{if $P\in\Har6(\RR^{24})$,}
          \end{cases}
      \end{equation*}
      where $c_2$ and $c_3$ are linearly independent linear forms on $\Har4(\RR^{24})$,
      and $c_4$ and $c_5$ are linearly independent linear forms on $\Har6(\RR^{24})$.
    \item If $\Lambda_2$ is empty, then
      \begin{equation*}
        \Theta_{\Lambda,P} =
          \begin{cases}
            \Th3^{24} - 48\,\Th3^{16}\De8 + 48\,\Th3^{8}\De8^2
                                                & \text{if $P=1$,} \\
            0                                   & \text{if $P\in\Har2(\RR^{24})$,} \\
            c_3(P)\,\Th3^{8}\De8^3              & \text{if $P\in\Har4(\RR^{24})$,} \\
            c_5(P)\,\Ph\Th3^{8}\De8^3           & \text{if $P\in\Har6(\RR^{24})$,}
          \end{cases}
      \end{equation*}
      where $c_3$ and $c_5$ are nonzero linear forms on $\Har4(\RR^{24})$ and $\Har6(\RR^{24})$ respectively.
  \end{enumerate}
\end{Lemma}

\begin{proof}[Idea of the proof]
  We use Proposition~\ref{CalcThetaSeries} for computing the theta series.
  The main difficulty is to show that, for example,
  in case~(ii),
  the coefficients $c_2$ and $c_3$ are linearly independant.
  (Similar cases are treated in the same way.)
  The main idea is to look at the theta series of the shadow (Proposition~\ref{ThetaSeriesOfShadow}):
  \begin{gather*}
      \Theta_{\Lambda,P} = c_2(P)\,\Th3^{16}\,\De8^2 + c_3(P)\,\Th3^{8}\De8^3 = c_2(P)\,q^2+O(q^3), \\
      \Theta_{\Sh(\Lambda),P} = c_2(P)\,\Th2^{16}\,\Sh\De8^2 + c_3(P)\,\Th2^{8}\Sh\De8^3
                         = -2^{4}\,c_3(P)\,q^2 + O(q^4).
  \end{gather*}
  So, we have
  \[
    c_2(P) = \sum_{x\in\Lambda_2} P(x)
    \quad\text{and}\quad 
    c_3(P) = -2^{-4}\sum_{x\in\Sh(\Lambda)_2} P(x).
  \]
  The shape of $\Lambda_2$ and $\Sh(\Lambda)_2$ is given
  by Proposition~\ref{OddSelfdualLatticesDimension24}.
  Finally, one can find two harmonic homogeneous polynomials of degree~$4$
  which give linearly independent values of $c_2(P)$ and $c_3(P)$.
\end{proof}

\medbreak

We can now conclude:
  
\begin{Theorem}\label{StrengthOddUnimodularRank24Lattices}
  Let $\Lambda$ be an odd selfdual lattice of rank~$24$ and of minimum~$2$.
  \begin{enumerate}
    \item If $\Lambda_2$ is strongly eutactic, then all nonnempty shells of $\Lambda$ and $\Sh(\Lambda)$
      are $3$-designs.
    \item If $\Lambda_2$ is not strongly eutactic, then the shells $\Lambda_m$ and $\Sh(\Lambda)_m$
      are $3$-designs for $m=4^a$, $a\ge1$.
  \end{enumerate}
\end{Theorem}

\begin{Remark}
  We have checked numerically that the shells of norm at most~$1200$
  of these lattices and of their shadows
  are not spherical designs of higher strength. 
\end{Remark}

\section{Other selfdual lattices up to rank~24}\label{SectOtherSelfdualLattices}

According to their classification
(\cite[Chap.\ 16~and~17]{ConwaySloane} and \cite{Bacher}),
the remaining selfdual lattices of rank at most~$24$ enter
in one of the three cases of the following Lemma:

\begin{Lemma}
  Let $\Lambda=\Z{p}\oplus L\sseq\RR^n$
  be a selfdual lattice with $\sigma(\Lambda)=n-\nobreak16$,
  where $L\sseq\RR^{N}$ is of minimum~$2$.
  \begin{enumerate}
    \item
      If $p=0$ and if $L_2$ is not strongly eutactic, then
      \begin{equation*}
        \Theta_{\Lambda,P} =
            c(P)\,\Ph\Th3^{n-16}\De8^2
           \qquad \text{if $P\in\Har2(\RR^n)$,}
      \end{equation*}
      where $c$ is a nonzero linear form on $\Har2(\RR^n)$.
    \item
      If $p\ge1$ and if $L_2$ is strongly eutactic of Coxeter number~$h$, then
      \begin{equation*}
        \Theta_{\Lambda,P} =
            c(P)\,\bigl(\Ph\Th3^{n-8}\De8 + (46-2N-h)\,\Ph\Th3^{n-16}\De8^2\bigr)
            \qquad \text{if $P\in\Har2(\RR^n)$,}
      \end{equation*}
      where $c$ is a nonzero linear form on $\Har2(\RR^n)$.
    \item
      If $p\ge1$ and if $L_2$ is not strongly eutactic, then
      \begin{equation*}
        \Theta_{\Lambda,P} =
            c_1(P)\,\Ph\Th3^{n-8}\De8 + c_2(P)\,\Ph\Th3^{n-16}\De8^2
            \qquad \text{if $P\in\Har2(\RR^n)$,}
      \end{equation*}
      where $c_1$ and $c_2$ are linearly independent linear forms on $\Har2(\RR^n)$.
  \end{enumerate}
\end{Lemma}

We have checked numerically that the shells of norm at most~$1200$
of these lattices and of their shadows
are not $3$-spherical designs.

\bigbreak
\section*{Appendix: The cubic lattices of rank 4 and~7}

The aim of this section is to prove the following result, without the use of modular forms:

\begin{Theorem}\label{TheoremCubicLattices4And7}\ \NOPAGEBREAK
  \begin{enumerate}
    \item
      Let $\Z{4}=\{x=(x_1,\ldots,x_4)\in\RR^4\mid x_i\in\ZZ\}$ be the cubic lattice of rank~$4$,
      and let $m$ be an even positive integer.
      Then the shell $(\Z{4})_m$ is a $5$-design.
    \item
      Let $\Z{7}=\{x=(x_1,\ldots,x_7)\in\RR^7\mid x_i\in\ZZ\}$ be the cubic lattice of rank~$7$,
      and let $m$ be a positive integer of the form $m=4^a(8b+3)$, $a,b\ge0$.
      Then the shell $(\Z{7})_m$ is a $5$-design.
  \end{enumerate}
\end{Theorem}

We begin with the proof of Claim~(i).
Let $C\sseq(\ZZ/2\ZZ)^4$ be the even weight code of length~$4$, which is defined by
\[
  c=(c_1,c_2,c_3,c_4) \in C  \iff c_1+c_2+c_3+c_4 = 0\in \ZZ/2\ZZ.
\]
For $y=(y_1,\dots,y_4)\in\Z{4}$ we write $\reduce{y}\in(\ZZ/2\ZZ)^4$ its class modulo $2\,\Z{4}$. 
Let $\Lambda$ be the sublattice of $\Z{4}$
consisting of elements $x\in\Z{4}$ such that $\reduce{x}\in C$.
It is a sublattice of index~$2$,
which is equivalent to the root lattice of~$\D_4$.

All shells of $\Lambda$ are $5$-designs.
This follows from the fact that
$\Lambda$ is invariant under the Weyl group $\Weyl(\F_4)$,
and there is no nonconstant harmonic polynomial of degree at most~$5$
that is invariant under the action of $\Weyl(\F_4)$
(\cite[Thm.~6.1]{GoethalsSeidel},
\cite[Thm.~3.12]{GoethalsSeidel81};
see also \cite[Sect.~4]{HarPac}).

Finally, for $m$ an even positive integer, we have $(\Z{4})_m=\Lambda_m$;
therefore $(\Z{4})_m$ is a spherical $5$-design.
This proves Claim~(i) of the Theorem.

\medbreak
Let us now show Claim~(ii).
Let $H\sseq(\ZZ/2\ZZ)^7$ be the Hamming code of length~$7$. Recall that it is a
linear code of minimal Hamming distance~$3$, containing
\begin{myitemize}
  \item one codeword $\mathbf{0}:=(0,0,0,0,0,0,0)$ of weight~$0$;
  \item seven codewords of weight~$3$;
  \item seven codewords of weight~$4$;
  \item one codeword $\mathbf{1}:=(1,1,1,1,1,1,1)$ of weight~$7$.
\end{myitemize}
The set $S$ of codewords of weight~$3$ forms a Steiner system $\mathrm{S}(2,3,7)$,
and the set of codewords of weight~$4$ is $S+\mathbf{1}$.
See for example \cite[Sect.~1.2, p.~7]{Ebeling}.

Let $\Lambda$ be the sublattice of $\Z{7}$
consisting of elements $x\in\Z{7}$ such that $\reduce{x}\in H$.
It is a sublattice of index~$8$, which is equivalent to $\sqrt{2}\,\E_7^\sharp$,
the rescaled weight lattice of $\E_7$.

As in the previous case, all shells of $\Lambda$ are $5$-designs,
because
$\Lambda$ is invariant under the Weyl group $\Weyl(\E_7)$,
and there is no harmonic polynomial of degree at most~$5$
that is invariant under the action of~$\Weyl(\E_7)$.

But contrarily to the previous case, no shell of $\Z{7}$
is equal to any shell of $\Lambda$.
Therefore, we need to look at the effect
of the action of some finite subgroup of
the orthogonal group of rank~$7$
on the shells of $\Lambda$.
  
We recall first what is a \emph{weighted} spherical design.

\begin{Definition}\label{DefWeightedDesign}
  A \emph{weighted spherical $t$-design} or \emph{spherical cubature formula of strength~$t$}
  is the data consisting of
  a nonempty finite subset~$X$ of $\Se_m$
  and a positive function $w:X\to\RR_{>0}$, $x\mapsto w_x$,
  such that $\sum_{x\in X}w_x=1$ and  such that
  the condition
  \begin{equation*}\tag{$C_{j}$}%
    \sum_{x\in X} w_x\,P(x)=0,\qquad \forall P\in\Har{j}(\RR^n)
  \end{equation*}
  holds for every integer~$j$ with $1\le j\le t$.
\end{Definition}

(For more on cubature formulae on spheres, see for example
\cite{Gaeta}.)

Occasionaly, we allow the weight function $w:X\to\RR$ to take the value zero.
A spherical design is a weighted spherical design with
constant weight function $w_x=1/\card{X}$.

\begin{Lemma}\label{LemmaWeightedDesign}
  Let $X\sseq\Se_m$ be a spherical $t$-design and let $G$ be a finite subgroup of
  the orthogonal group $\Orth(n)$.
  Then $GX$ is a weighted spherical $t$-design for the weight function
  \[
    w_y=\frac{\Card{\{\sigma\in G \mid \sigma y\in X\}}}{\card{G}\,\card{X}}, \qquad y\in GX.
  \]
\end{Lemma}

\begin{proof}
  For $P\in\Har{j}(\RR^n)$, $1\le j\le t$, we have
  \[
    \sum_{y\in GX} w_y P(y)
    = \frac{1}{\card{G}\card{X}}\sum_{\substack{x\in X\\ \sigma\in G}} P(\sigma x) 
    = \frac{1}{\card{G}\card{X}}\sum_{\sigma\in G} \biggl(\sum_{x\in X} (P\circ \sigma)(x)\biggr).
  \]
  The last term is zero, because $P\circ\sigma\in\Har{j}(\RR^n)$
  and $X$ is a spherical $t$-design.
\end{proof}

We apply the previous Lemma to the subgroup $G=\Aut(\Z{7})\simeq(\ZZ/2\ZZ)^7\rtimes\mathfrak{S}_7$ of $\Orth(7)$,
which consists of transformations of the form
\[
  (x_1,\dots,x_7)\mapsto  (\epsilon_1 x_{\pi(1)},\dots,\epsilon_7 x_{\pi(7)})
\]
where $\epsilon_i\in\{\pm1\}$ and $\pi$ is a permutation of the set $\{1,2,\dots,7\}$.

Let $c\in\{0,1,\dots,7\}$, and let $m$ be a positive integer such that
$m\equiv c\bmod 8$.
For $y=(y_1,\dots,y_7)\in(\Z{7})_m$,
let $W(y)$ be the weight of $\reduce{y}\in(\ZZ/2\ZZ)^7$,
that is the number of coordinates~$i$ such that $y_i$ is odd.
(Caution: the word ``weight'' has two different meanings in this Appendix.)
The condition $y_1^2+\dots+y_7^2\equiv c\bmod 8$
implies that $W(y)$ takes the values indicated in the following table:
$$
  \newcommand*\cm{\mathbin{,}}
  \vcenter{\ialign{%
    \strut\hfil\ $#$\ \hfil&\vrule\hfil\ $#$\ \hfil\cr
    c & W(y) \cr
    \noalign{\hrule}
    0 & 0\cm4 \cr
    1 & 1\cm5 \cr
    2 & 2\cm6 \cr
    3 & 3   \cr
  }}
  \qquad\qquad
  \vcenter{\ialign{%
    \strut\hfil\ $#$\ \hfil&\vrule\hfil\ $#$\ \hfil\cr
    c & W(y) \cr
    \noalign{\hrule}
    4 & 0\cm4 \cr
    5 & 1\cm5 \cr
    6 & 2\cm6 \cr
    7 & 3\cm7 \cr
  }}
$$
Now, the quantity
\[
  \lambda(y):=\frac{\Card{\{\sigma\in G \mid \sigma y\in \Lambda_m\}}}{\card{G}},
  \qquad\text{for $y\in(\ZZ^n)_m$}
\]
depends only of $W(y)$; it
is given by the following table:
$$
  \vcenter{\ialign{%
    \strut\hfil\ $#$\ \hfil&\vrule\hfil\ $#$\ \hfil\cr
    W(y) & \lambda(y) \cr
    \noalign{\hrule}
    0 & 1  \cr
    1 & 0  \cr
    2 & 0  \cr
    3 & 1/5\cr
  }}
  \qquad\qquad
  \vcenter{\ialign{%
    \strut\hfil\ $#$\ \hfil&\vrule\hfil\ $#$\ \hfil\cr
    W(y) & \lambda(y) \cr
    \noalign{\hrule}
    4 & 1/5 \cr
    5 & 0   \cr
    6 & 0   \cr
    7 & 1   \cr
  }}
$$
Applying Lemma~\ref{LemmaWeightedDesign}, we find that
\begin{myitemize}
  \item for $m\equiv0\bmod8$ and $m\equiv4\bmod8$,
    the shell $(\Z{7})_m$ is a weighted spherical $5$-design for the weight function
    \[
      w_y= \frac{1}{\card{\Lambda_m}}
      \text{ if $W(y)=0$,  and  }
      w_y= \frac{1}{5\card{\Lambda_m}}
      \text{ if $W(y)=4$;} 
    \]
  \item for $m\equiv3\bmod8$, the shell $(\Z{7})_m$ is a spherical $5$-design;
  \item for $m\equiv7\bmod8$, the shell $(\Z{7})_m$ is a weighted spherical $5$-design for the weight function
    \[
      w_y= \frac{1}{\card{\Lambda_m}}
      \text{ if $W(y)=7$,  and  }
      w_y= \frac{1}{5\card{\Lambda_m}}
      \text{ if $W(y)=3$.} 
    \]
\end{myitemize}

In order to achieve the proof of Theorem~\ref{TheoremCubicLattices4And7}, Claim~(ii),
it remains to show the following statement:
\begin{quotation}
  \noindent Let $m$ be a positive integer.
  If $(\Z{7})_m$ is a $5$-design, then $(\Z{7})_{4m}$ is also a $5$-design.
\end{quotation}
The proof is the following.
We write $(\Z{7})_{4m}$ as
\[
  (\Z{7})_{4m} = 2\,(\Z{7})_m \sqcup Q
\]
where
$2\,(\Z{7})_m$ is the shell of norm~$m$ rescaled by a factor~$2$,
and $Q$ contains the elements $y\in(\Z{7})_{4m}$
such that $W(y)=4$.
We have shown that $(\Z{7})_{4m}$
is a weighted $5$-design for the weight function
\[
      w_y = \frac{1}{\card{\Lambda_m}}
      \text{ if $y\in2\,(\Z{7})_m$,  and  }
      w_y = \frac{1}{5\card{\Lambda_m}}
      \text{ if $y\in Q$.} 
\]
But since $(\Z{7})_m$ is a $5$-design by hypothesis,
then $(\Z{7})_{4m}$
is a weighted $5$-design for the weight function
\[
      \widetilde{w}_y = \frac{1}{\card{(\Z{7})_m}}
      \text{ if $y\in2\,(\Z{7})_m=0$,  and  }
      \widetilde{w}_y = 0 \text{ if $y\in Q$.} 
\]

It is evident that, if $X$ is a weighted $t$-design for
two weight functions $w$ and $\widetilde{w}$,
then it is a weighted $t$-design for every convex linear combination
of $w$ and $\widetilde{w}$.

In particular, since, in our case, there is a suitable
convex linear combination
of $w$ and $\widetilde{w}$ which is constant on $(\Z{7})_{4m}$,
the shell $(\Z{7})_{4m}$ is indeed a $5$-design.

\smallskip

Remark that our proof also shows that the shells of $\Z{7}$ of
norm $m\equiv 7\bmod8$ and those of norm $m\equiv0\bmod4$
are weighted spherical $5$-designs,
although they are not spherical $5$-designs in general.

  
\bigskip
\noindent\textbf{Acknowledgements.}
The author thanks Boris \textsc{Venkov},
who was at the origin of this work,
and who introduced the author to the use of modular forms
in connection with lattices.
Pierre \textsc{de la Harpe} also contributed to this work.

The author acknowledges support from the \emph{Swiss National Science Foundation}.



\begin{thebibliography}{HarPac04b}
  

  \bibitem[Bach97]{Bacher}
    R.~Bacher,
    \emph{Tables de r\'eseaux entiers construits comme $k$-voisins de $\ZZ^n$},
    J.\ Th\'eor.\ Nombres Bordeaux \textbf{9} (2) (1997) 479--497.

  \bibitem[BacVen01]{BachocVenkov}
    Ch.~Bachoc, B.B.~Venkov,
    \emph{Modular forms, lattices and spherical designs},
    chap.~2 of~\cite{MartinetVenkov}

  \bibitem[BekHar02]{BekkaHarpe}
    B.~Bekka, P.~de la Harpe,
    \emph{Irreducibility of unitary group representations and reproducing kernels Hilbert spaces},
    Expo.\ Math.\ \textbf{21} (2) (2003) 115--149.
        
  \bibitem[Borc84]{BorcherdsOdd24Ori}
    R.E. Borcherds,
    \emph{The Leech lattice and other latices},
    Ph.D.\ dissertation,
    Univ.\ of Cambridge 1984.

  \bibitem[BorcCS]{BorcherdsOdd24}
    R.E. Borcherds,
    \emph{The 24-dimensional odd unimodular lattices},
    Chap.~17 of~\cite{ConwaySloane}.
  
  \bibitem[Bour81]{BourbakiLie}
      N.~Bourbaki,
      \emph{Groupes et alg\`ebres de Lie},
      Masson (Paris) 1981.

  \bibitem[ConSlo99]{ConwaySloane}
      J.H.~Conway, N.J.A.~Sloane,
      \emph{Sphere Packings, Lattices and Groups},
      third edition, Springer (New York) 1999.

  \bibitem[DeGoSe77]{DelsarteGoethalsSeidel}
    P.~Delsarte, J.-M.~Goethals, J.J.~Seidel,
    \emph{Spherical codes and designs},
    Geometriae Dedicata  \textbf{6} (1977) 363--388.
    
  \bibitem[Ebel94]{Ebeling}
      W.~Ebeling,
      \emph{Lattices and codes, a course partially based on lectures by F.~Hirzenbruch},
      Vieweg (Braunschweig) 1994.
      Second revised edition: 2002.
      
  \bibitem[Elki95a]{Elkies0}
      N.D.~Elkies,
      \emph{A characterization of the $\Z{n}$-lattice},
      Math.\ Res.\ Lett.\ \textbf{2} (3) (1995) 321--326.

  \bibitem[Elki95b]{Elkies8}
      N.D.~Elkies,
      \emph{Lattices and codes with long shadows},
      Math.\ Res.\ Lett.\ \textbf{2} (5) (1995) 643--645.


  \bibitem[GoeSei79]{GoethalsSeidel}
    J.-M.~Goethals, J.J.~Seidel,
    \emph{Spherical Designs},
    Proc.\ Sympos.\ Pure Math., XXXIV, 255--272,
    Amer.\ Math.\ Soc. (Providence R.I.) 1979.
    
  \bibitem[GoeSei81]{GoethalsSeidel81}
    J.-M.~Goethals, J.J.~Seidel,
    \emph{Cubature formulae, polytopes, and spherical designs},
    in: \emph{The geometric vein, the Coxeter festschrift},
    Springer (New York, Berlin) (1981) 203--218.
    
  \bibitem[HarPac04a]{HarPac}
    P.~de la Harpe, C.~Pache,
    \emph{Spherical designs and finite group representations
          (some results of E.~Bannai)},
    European J.\ Combin.\  \textbf{25} (2) (2004) 213--227.

  \bibitem[HarPac04b]{Gaeta}
    P.~de la Harpe, C.~Pache,
    \emph{Cubature formulas, geometrical designs, reproducing kernels, and Markov operators},
    in preparation.

  \bibitem[MartV01]{MartinetVenkov}
      J.~Martinet ed.,
      \emph{R\'eseaux euclidiens, designs sph\'eriques et formes modulaires,
          Autour des travaux de B.~Venkov}, 
      L'Enseignement Math\'ematique,
      monographie n$\rm^o$~37
      (Gen\`eve) 2001.

  \bibitem[NebVen00]{NebeVenkov10}
    G.~Nebe, B.B.~Venkov,
    \emph{The strongly perfect lattices of dimension 10},
    J.\ Th\'eor.\ Nombres Bordeaux \textbf{12} (2) (2000) 503--518.
  
 \bibitem[Rank77]{Rankin}
      R.A.~Rankin,
      \emph{Modular forms and functions},
      Cambridge University Press, 1977.
      
  \bibitem[SchSch99]{ScharlauSchulze-Pillot}
      R.~Scharlau, R.~Schulze-Pillot, 
      \emph{Extremal lattices},
      in: \emph{Algorithmic algebra and number theory (Heidelberg, 1997)},
      Springer (Berlin) (1999) 139--170.

  \bibitem[Serr85]{Serre85}
      J.-P.~Serre,
      \emph{Sur la lacunarit\'e des puissances de $\eta$},
      Glasgow Math.\ J.\ \textbf{27} (1985) 203--221
      [=~\OE uvres, Volume~IV, 66--84, see also 640].

  \bibitem[Venk78]{Venkov24}
      B.B.~Venkov
      \emph{The classification of integral even unimodular 24-\hskip0pt\relax dimensional quadratic forms},
      Trudy Mat.\ Inst.\ Steklova \textbf{148} (1978) 65--76
      =~Proc.\ Steklov Inst.\ Math.\ n$\rm^o$4 (1980) 63--74 \\
      $\simeq$~\emph{Even Unimodular 24-Dimensional Lattices},
      Chap.~18 of \cite{ConwaySloane}.

  \bibitem[Venk84]{VenkovExtremal}
      B.B.~Venkov,
      \emph{Even unimodular extremal lattices},
      Tr.\ Mat.\ Inst.\ Steklova \textbf{165} (1984) 43--48 (Russian)
      = Proc.\ Steklov Inst.\ Math.\ \textbf{165} (3) (1985) 47--52.
   
  \bibitem[VenMar01]{VenkovMartinet1}
      B.B.~Venkov (notes by J.~Martinet),
      \emph{R\'eseaux et designs sph\'eriques},
      Chap~1 of \cite{MartinetVenkov}.

  \bibitem[Vile68]{Vilenkin}
      N.Ya.~Vilenkin,      
      \emph{\textcyrillic{Spetsialp1nye funktsii i teoriya predstavleni\u i grupp}},       
      second edition, Nauka (Moscow) 1991 \\
      =~\emph{Special Functions and the Theory of Group Representations},
      Transl.\ Math.\ Monographs 22,
      Amer.\ Math.\ Soc. (Providence R.I.) 1968 \\
      =~\emph{Fonctions sp\'eciales et th\'eorie de la repr\'esentation des groupes},
      Dunod (Paris) 1969.
      
  \bibitem[Voro08]{Voronoi}
      G.~Voronoi,
      \emph{Nouvelles applications des param\`etres continus \`a la th\'eorie
      des formes quadratiques~:
      1.\ Sur quelques propri\'et\'es des formes quadratiques positives parfaites},
      J.\ reine angew.\ Math.\ \textbf{133} (1908) 97--178.
       
\end{thebibliography}
\end{document}